\documentclass[twocolumn]{autart}    

\usepackage{graphics} 
\usepackage{epsfig} 
\usepackage{mathptmx} 
\usepackage{times} 
\usepackage{amsmath} 
\usepackage{amssymb}  
\usepackage{amsbsy,bm,fixmath}

\usepackage{graphicx}          

\newtheorem{theorem}{Theorem}
\newtheorem{definition}[theorem]{Definition}
\newtheorem{lemma}[theorem]{Lemma}
\newtheorem{proposition}[theorem]{Proposition}

\newtheorem{example}{Example}
\newtheorem{remark}{Remark}

\def\mbA{\mathbold{A}}
\def\mbB{\mathbold{B}}
\def\mbD{\mathbold{D}}
\def\mbP{\mathbold{P}}
\def\mbx{\mathbold{x}}
\def\mbQ{\mathbold{Q}}
\def\mbS{\mathbold{S}}
\def\mbK{\mathbold{K}}

\def\mbr{\mathbold{r}}

\def\proof{{\it Proof}.\quad}
\def\endproof{\hfill$\Box$}

\begin{document}

\begin{frontmatter}

\title{Linear quadratic mean field games with a major player:\\ The multi-scale approach \thanksref{footnoteinfo}} 

\thanks[footnoteinfo]{This paper was not presented at any IFAC
meeting. 
 This work was supported by
 the National Science Foundation of China (No.11601489), Startup Research Fund of Zhengzhou University (No.129-51090091), Outstanding Young Talent Research Fund of Zhengzhou University (No.129-32210453), Natural Sciences and Engineering Research Council (NSERC) of Canada. Submitted to Automatica, Jan 2019; revised Aug 2019. This version contains a more detailed Sec. 5 than the revised journal submission. Corresponding author: M. Huang.}

\author[Ma]{Yan Ma}\ead{mayan203@zzu.edu.cn},    
\author[Huang]{Minyi Huang}\ead{mhuang@math.carleton.ca}  

\address[Ma]{School of
Mathematics and Statistics, Zhengzhou University, Zhengzhou, 450001, China}  
\address[Huang]{School of
Mathematics and Statistics, Carleton University, Ottawa, ON K1S 5B6,
Canada}             


\begin{abstract}                          
This paper considers linear quadratic  (LQ) mean field games with a major player and analyzes an asymptotic solvability problem. It starts with a large-scale system of coupled dynamic programming equations and applies a re-scaling technique introduced in Huang and Zhou (2018a, 2018b) to derive a set of Riccati equations in lower dimensions, the solvability of which determines the necessary and sufficient condition for asymptotic solvability. We next derive the mean field limit of the strategies and the value functions. Finally, we show that the two decentralized strategies    can be interpreted as the best responses of a major player and a representative minor player embedded in an infinite population, which have the property of consistent mean field approximations.
\end{abstract}

\end{frontmatter}

\section{Introduction}

Mean field game theory has undergone a phenomenal growth. It provides a powerful methodology
for handling complexity in noncooperative mean field decision problems.
The readers are referred to (Caines, Huang, and Malham\'e, 2017) for an overview.
 Most  existing analysis has been developed based on two routes called the direct approach and the fixed point approach. By the direct approach, one starts by formally solving an $N$-player game to obtain a large coupled solution equation system, and  next  derives a simple limiting equation system by taking $N\to \infty$; see (Lasry and Lions, 2007) for the limit consisting of a Hamilton-Jacobi-Bellman (HJB) equation and a Fokker-Planck-Kolmogorov (FPK) equation. By the fixed point approach, one determines the best response of a representative agent to a mean field of an infinite population, and next all the agents' best responses should regenerate that mean field (Huang, Malham\'e, and Caines, 2006).
This procedure formalizes a fixed point problem, which can be solved and further used to design decentralized strategies. For LQ mean field games, the recent work (Huang and Zhou, 2018b) shows the exact relationship of the two approaches. In general, the fixed point approach has more flexibilities and can be implemented in diverse models (Huang, Caines, and Malham\'e, 2007; Li and Zhang, 2008; Bensoussan et al, 2013; Huang and Ma, 2016; Carmona and Delarue, 2018).  Further convergence analysis in the direct approach can be found in  (Cardaliaguet et al, 2015; Lacker, 2016; Fischer, 2017). Mean field games have found applications in traffic routing (Bauso, Zhang, and  Papachristodoulou, 2017), smart grids (Couillet, et al, 2012; Ma, Callaway, and Hiskens, 2013;
Kizilkale, Salhab, and Malham\'e, 2019) and production planning (Wang and Huang, 2019), among others.  A notable feature of the early literature of mean
field games is that all players in the model are comparably small, and can be called peers.

Huang (2010)  introduces an LQ mean field game model with a major player which has strong influence. A motivating example is the interaction between a large corporation and many much smaller firms. There has been a rapid increase of literature on mean field games with a major and many minor players. In the setting of LQ models, Nguyen and Huang (2012a) consider continuum parametrized minor players, and Nguyen and Huang (2012b) extend to mass behavior directly impacted by the major player.  Kordonis and Papavassilopoulos  (2015) analyze minor players with random entrance. Major players with leadership are studied by Bensoussan et al (2017), Moon and Basar (2018). Partial state observation is considered by Caines and Kizilkale (2017),  Firoozi and Caines (2015).  Huang, Wang and Wu (2016) take linear backward stochastic differential equations to model the dynamics of the players.
Huang, Jaimungal, and Nourian (2015)  present an application of
the major player mean field game theory to an optimal execution model with an institutional trader and a large number of small traders.

 Major-minor player games with nonlinear diffusion dynamics are an important class of modelling; see Nourian and Caines (2013), Buckdahn, Li and Peng (2014), Bensoussan, Chau and Yam (2016),  Carmona and Zhu (2016). Leader-follower interaction is adopted by Bensoussan et al (2015), Fu and Horst (2018). To deal with this nonlinear modelling, forward-backward stochastic differential equations provide a vital analytical tool. Sen and Caines (2016) apply nonlinear filtering when the major player's state is partially observed. More recently, Lasry and Lions (2018) introduce master equations for mean field games with major and minor players. They may be viewed as a pair of abstract dynamic programming equations. Cardadiaguet, Cirant, and Porretta (2018) prove the convergence of the Nash equilibria by use of the master equations when the number of minor players tends to infinity.
A mean field principal-agent model is  formulated by  Elie,  Mastrolia,   and Possamai (2019).
For major player models with discrete states, see
(Huang 2012; Carmona and Wang, 2017; Kolokoltsov, 2017).

 Huang (2010) applies a state space augmentation approach by adding the mean field dynamics into the two decision problems, one for the major player and one for a representative minor player. This Markovianizes the problem and enables the use of dynamic programming. The procedure of Huang (2010) is based on  the fixed point approach and the associated consistent mean field approximations, and that work only assumes
existence of the solution.

This paper analyzes the LQ mean field game with a major player and homogeneous (or symmetric) minor players and takes the direct approach by starting with the solution  for $N+1$ players.
 Specifically, we will extend an asymptotic solvability notion introduced in a recent work Huang and Zhou (2018a) for LQ mean field games without a major player. With or without a major player,  asymptotic solvability can be informally stated as the existence of Nash equilibria with complete state information for all sufficiently large population sizes, in addition to some boundedness property of the solution.
We exploit the multi-scale nature of the optimization problem and  use a re-scaling method in Huang and Zhou (2018a, 2018b) so
that the key
information in some higher order terms, as  components  in the
solution matrices of $N+1$ coupled Riccati equations, can be captured. We  derive the necessary and sufficient condition for asymptotic solvability and evaluate the value function. The re-scaling method gives a set of ordinary differential equations  (ODEs) for nine matrix functions. To reveal the special structure underlying these functions, we will further
relate them to the best responses of the major player and a representative minor player staying an infinite population, where consistent mean field
approximations hold. The latter is a key feature of the fixed point approach in mean field games.  Our mean field limit analysis shares similarity to (Cardadiaguet, Cirant, and Porretta, 2018) which performs convergence analysis in a nonlinear system via the master equation. But we explicitly exploit the multi-scale phenomena in our model to identify a lower dimensional object which  governs the asymptotic behavior of the system when the number of minor players tends to infinity. Similar methods  appear in the statistical physics literature on mean field models (Ott and Antonsen, 2008; Pazo and Montbrio, 2014).

 We mention other related LQ models of finding mean field limits via analyzing large scale equations. Papavassilopoulos (2014)  uses  large algebraic Riccati equations in mean field games and analyzes existence by an implicit function theorem. Priuli (2015) considers coupled HJB and FPK equations with decentralized information.
 Mean field social optimal control is analyzed in (Huang 2003, Chap. 6; Herty, Pareschi, and Steffensen, 2015) via large Riccati equations.

The organization of the paper is as follows.
Section \ref{sec:mod} describes the LQ Nash game with $N+1$ players together with its solution via dynamic programming and Riccati equations.
Section \ref{sec:asy} extends the formulation of  asymptotic  solvability
in Huang and Zhou (2018a, 2018b) to the LQ model with a major player. Section \ref{sec:cos} presents further mean field limits
and the performance. Section \ref{sec:two} formulates two optimal control problems under a mean field generated by an infinite number of minor players and addresses the relation to the asymptotic solvability problem. Numerical examples are presented in Section \ref{sec:num}.
  Section \ref{sec:con} concludes the paper.

{\it Notation:} For symmetric matrix $S\ge 0$, we may write $x^TSx=|x|_S^2$.  For a matrix $ Z=(z_{jk}) \in \mathbb{R}^{l\times m}$, denote the $l_1$-norm $\| Z\|_{l_1}=\sum_{j,k}|z_{jk}|$. Let the function $g(\delta, x)$ be defined for $x$ in a subset $D_g$ of a
Euclidean space and parameter $\delta\in (0,p]$ for some $p>0$. We say $g$ is compactly of $O(\delta)$ if for each compact subset $D_0\subset D_g$, there exists a constant $c_0$ depending on $D_0$ such that $\sup_{x\in D_0} |g(\delta,x)|\le c_0 \delta$.

\section{The LQ game with major and minor players}
\label{sec:mod}

We consider the LQ game with a major player ${\mathcal A}_0$
and $N$
minor players ${\mathcal A}_i$,  $1\le i\le N$.
At time $t\geq 0$, the states of ${\mathcal A}_0$ and ${\mathcal A}_i$ are, respectively,
denoted by $X_0(t)$ and $X_i(t)$,  $1\le i\le N$.
 The dynamics of the $N +1$ players are given by a system of linear stochastic
differential equations (SDEs):
\begin{align}
dX_0(t)=\  &\big(A_0X_0(t)+B_0u_0(t)+F_0X^{(N)}(t)\big)dt \nonumber \\
& +D_0dW_0(t),\label{stateX0}\\
dX_i(t)=\  &\big(AX_i(t)+Bu_i(t)+FX^{(N)}(t)+GX_0(t)\big)dt\nonumber \\
& +DdW_i(t), \quad 1\le i\le N,\quad t\ge 0,\label{stateXi}
\end{align}
where we have   state  $X_i\in \mathbb{R}^n$, control
$u_i\in\mathbb{R}^{n_1}$, and
 $ X^{(N)}=\frac{1}{N}\sum_{k=1}^N X_k(t)$.
 The initial states $\{X_j(0), 0\le j\le N\}$ are independent with $EX_j(0)=x_j(0)$ and finite second moment. The $N+1$ standard  $n_2$-dimensional Brownian motions $\{W_j, 0\le j\le N\}$ are independent and also independent of the  initial states. The deterministic constant matrices
$A_0$, $A$, $B_0$, $B$, $D_0$, $D$, $F_0$, $F$, $G$ have compatible
dimensions. Denote $u=(u_0, \cdots, u_N)$.
The costs of players ${\mathcal A}_k$, $0\le k\le N$,  are given by
\begin{align}
J_0(u) =\  & E\int_0^T \left[  |X_0(t)-\Gamma_0 X^{(N)}(t)-\eta_0|_{Q_0}^2+|u_0(t)|_{R_0}^2\right]dt \nonumber\\
&+E|X_0(T)-\Gamma_{0f} X^{(N)}(T)-\eta_{0f}|_{Q_{0f}}^2,\label{costJ0} \\
J_i(u) =\ &  E\int_0^T \left[  |X_i(t)-\Gamma_1 X_0(t)-\Gamma_2 X^{(N)}(t)-\eta|_{Q}^2+|u_i(t)|^2_R \right]dt\nonumber\\
&\hskip -0.5cm+E|X_i(T)-\Gamma_{1f} X_0(T)-\Gamma_{2f} X^{(N)}(T)-\eta_f|_{Q_f}^2, \label{costJi} \\
& 1\le i\le N. \nonumber
\end{align}
The deterministic constant matrices (or vectors)
 $Q_0$, $\Gamma_0$, $\eta_0$, $R_0$, $Q_{0f}$, $\Gamma_{0f}$, $\eta_{0f}$, $Q$, $\Gamma_1$, $\Gamma_2$, $\eta$, $R$, $Q_{f}$, $\Gamma_{1f}$, $\Gamma_{2f}$, $\eta_f$ above
have compatible dimensions, and $Q_{0}\geq 0$, $Q_{0f}\geq 0$, $Q\geq 0$, $Q_f\geq 0$, $R_0>0$, $R>0$.
For notational simplicity, we only consider constant parameters. Our analysis can be easily extended to
the case of time-dependent parameters.
Define
\begin{align}
&\mbx=(\mbx_0^T,\mbx_1^T,\cdots,\mbx_N^T)^T\in\mathbb{R}^{(N+1)n},
\nonumber\\
&X(t)= \begin{bmatrix} X_0(t) \\ \vdots \\ X_N(t)
\end{bmatrix}\in \mathbb{R}^{(N+1)n},\nonumber
\quad W(t)=\begin{bmatrix}W_0(t) \\ \vdots \\ W_N(t)
\end{bmatrix}\in \mathbb{R}^{(N+1)n_2},\nonumber   \\
&\widehat \mbA = \mbox{diag}[A_0, A, \cdots, A]+\left[
\begin{matrix}
& 0,&{\bf 1}_{1\times n}\otimes \frac{F_0}{N}\\
& {\bf 1}_{n\times 1}\otimes G, &{\bf 1}_{n\times n}\otimes
\frac{F}{N}
\end{matrix}\right],\nonumber\\
&\widehat \mbD =\mbox{diag}[D_0, D, \cdots, D]\in \mathbb{R}^{(N+1)n\times (N+1)n_2},\nonumber \\
&\mbB_0 = e_1^{N+1}\otimes B_0\in\mathbb{R}^{(N+1)n\times n_1}, \nonumber\\
&\mbB_{k} = e_{k+1}^{N+1}\otimes B\in\mathbb{R}^{(N+1)n\times n_1},\qquad 1 \leq k \leq N. \nonumber
\end{align}
We denote by ${\bf 1}_{k\times l}$  a $k\times l$ matrix with all entries equal to 1, and by the column  vectors  $\{e_1^k, \ldots, e_k^k\}$  the canonical basis of $\mathbb {R}^k$. For instance, $e_1^k=[1, 0 ,\cdots, 0]^T\in \mathbb{R}^k$.
 For matrices $K=(k_{ij})\in \mathbb{R}^{l_1\times l_2}$, $\hat K\in \mathbb{R}^{l_3 \times l_4}$, the Kronecker product
 $K\otimes \hat K= (k_{ij} \hat K)_{1\le i\le l_1,1\le j\le l_2}\in
 \mathbb{R}^{(l_1l_3)\times (l_2 l_4)}$.
   We may use a subscript $n$ to indicate the identity matrix $I_n$ to be $n\times n$.

Now we write \eqref{stateX0} and \eqref{stateXi} in the form
\begin{align}\label{bigx}
dX(s)=\Big(\widehat \mbA X(s)+\sum_{k=0}^N \mbB_k u_k(s)\Big)dt+\widehat \mbD dW(s),
\end{align}
where $s\ge 0$.
We consider closed-loop perfect state (CLPS) information so that
$X(s)$ is observed by each player, and look for Nash strategies in this section. Let $u_{-k}$ denote the strategies of all players other than ${\mathcal A}_k$.
 A set of strategies
 $(\hat u_0,\cdots, \hat u_N)$ is a Nash equilibrium
if for any $0\le k\le N$, we have
$$
J_k(\hat u_k, \hat u_{-k})\le J_k(u_k, \hat u_{-k}),
$$
for any state feedback based strategy $u_k$ which together with $\hat u_{-k}$ ensures a unique solution of $X(s)$ on $[0,T]$.
Denote the value function of ${\mathcal A}_j$  by $V_j(t,\mbx)$, $0\le j\le N$, which corresponds to the initial time-state pair $(t,\mbx)$ in \eqref{bigx}, i.e., $X(t)= \mbx$ at the initial time $t$, and can be interpreted as $J_j$ evaluated on the time interval $[t,T]$ under the set of Nash strategies.
The set of value functions is determined by the system of  HJB equations
\begin{align}\label{DE0}
&0= \frac{\partial V_0}{\partial t}+\min_{u_0\in {\mathbb R}^{n_1}}\Big(\frac{\partial^T V_0}{\partial \mbx }\Big(\widehat {\mbA}\mbx+\sum_{k=0}^N {\mbB}_k u_k\Big)
+u_0^T R_0 u_0 \nonumber \\
   &\qquad+|\mbx_0-\Gamma_{0} \mbx^{(N)}-\eta_{0}|_{Q_{0}}^2  +\frac{1}{2}\mbox{Tr}\big
   ({\widehat{\mbD}^T (V_0)_{\mbx\mbx} \widehat{\mbD}}\big)\Big), \\
&V_0(T,\mbx)=|\mbx_0-\Gamma_{0f} \mbx^{(N)}-\eta_{0f}|_{Q_{0f}}^2, \nonumber
\end{align}
and
\begin{align}\label{DE}
&0= \frac{\partial V_i}{\partial t}+\min_{u_i\in {\mathbb R}^{n_1}}\Big(\frac{\partial^T V_i}{\partial \mbx }\Big(\widehat {\mbA}\mbx+\sum_{k=0}^N {\mbB}_k u_k\Big)
+u_i^T R u_i \nonumber \\
   &\qquad+|\mbx_i-\Gamma_1 \mbx_0-\Gamma_2 \mbx^{(N)}-\eta|_{Q}^2  +\frac{1}{2}\mbox{Tr}\big
   ({\widehat{\mbD}^T (V_i)_{\mbx\mbx} \widehat{\mbD}}\big)\Big), \\
&V_i(T,\mbx)=|\mbx_i-\Gamma_{1f} \mbx_0-\Gamma_{2f} \mbx^{(N)}-\eta_f|_{Q_f}^2, \quad  1\le i\le N, \nonumber
\end{align}
where $\mbx^{(N)}=(1/N)\sum_{i=1}^N \mbx_i$
and the minimizers are
\begin{align}
&u_0=-\frac{1}{2} R_0^{-1} {\mbB}_0^T \frac{\partial V_0}
 {\partial \mbx},\\
 &u_i=-\frac{1}{2} R^{-1} {\mbB}_i^T \frac{\partial V_i}
 {\partial \mbx},\quad 1\leq i \leq N.
\end{align}
 Next we substitute $u_0$ and $u_i$ into \eqref{DE0} and \eqref{DE}:
 \begin{align}\label{DE2}
0 = &\frac{\partial V_0}{\partial t}+\frac{\partial^T V_0}{\partial \mbx}\Big({\widehat \mbA}\mbx
-\frac{1}{2}{\mbB}_0 R_0^{-1} {\mbB}_0^T \frac{\partial V_0}{\partial \mbx}\nonumber \\
& -\frac{1}{2}\sum_{k=1}^N {\mbB}_k R^{-1} {\mbB}_k^T \frac{\partial V_k}{\partial \mbx}\Big)
+|\mbx_0-\Gamma_{0} \mbx^{(N)}-\eta_{0}|_{Q_{0}}^2   \nonumber \\
   &+\frac{1}{4}\frac{\partial^T V_0}{\partial \mbx}{\mbB}_0 R_0^{-1} \mbB_0^T\frac{\partial V_0}{\partial \mbx}+\frac{1}{2}\mbox{Tr}\big({{\widehat \mbD}^T (V_0)_{\mbx\mbx} {\widehat \mbD}}\big),
\end{align}
and
\begin{align}\label{DE3}
0 = &\frac{\partial V_i}{\partial t}+\frac{\partial^T V_i}{\partial \mbx}\Big({\widehat \mbA}\mbx
-\frac{1}{2}{\mbB}_0 R_0^{-1} {\mbB}_0^T \frac{\partial V_0}{\partial \mbx}\nonumber \\
&-\frac{1}{2}\sum_{k=1}^N {\mbB}_k R^{-1} {\mbB}_k^T \frac{\partial V_k}{\partial \mbx}\Big)
+|\mbx_i-\Gamma_1 \mbx_0-\Gamma_2 \mbx^{(N)}-\eta|_{Q}^2 \nonumber \\
   &+\frac{1}{4}\frac{\partial^T V_i}{\partial \mbx}{\mbB}_i R^{-1} \mbB_i^T\frac{\partial V_i}{\partial \mbx} \nonumber \\
   &+\frac{1}{2}\mbox{Tr}\big({{\widehat \mbD}^T (V_i)_{\mbx\mbx} {\widehat \mbD}}\big),\qquad 1\leq i \leq N.
\end{align}

Suppose $V_{j}(t,\mbx),~0\leq j\leq N,$ has the following form
\begin{align}\label{Vform}
&V_{j}(t,\mbx)=\mbx^T {\mbP}_{j}(t) \mbx+2{\mbS}_{j}^T(t) \mbx+\mbr_{j}(t),
\end{align}
where $\mbP_{j}(t)$ is symmetric.
Then
\begin{align}  \label{1st}
\frac{\partial V_{j}}{\partial \mbx}=2{\mbP}_{j}(t)\mbx+2{\mbS}_{j}(t),\quad
\frac{\partial^2 V_{j}}{\partial \mbx^2}=2\mbP_{j}(t).
\end{align}

Denote
\begin{align}
&\mbK_{0}= [{I_n},0,\cdots,0]-\tfrac{1}{N}
[0,\Gamma_{0},\cdots,\Gamma_{0}],\nonumber  \\
& \mbQ_{0}=\mbK_{0}^T Q_{0} \mbK_{0},\nonumber  \\
&\mbK_{i}=[0,0,\cdots,{I}_{n},0,\cdots,0]-[\Gamma_{1},0,\cdots,0]\nonumber\\
&\quad\quad\quad -\tfrac{1}{N}[0,\Gamma_{2},\cdots,\Gamma_{2}], \label{kig}\\
& \mbQ_{i}=\mbK_{i}^T Q \mbK_{i}, \nonumber
\end{align}
where $I_n$ is the $(i+1)$th submatrix in \eqref{kig}. We have $\mbK_0,\mbK_i\in \mathbb {R}^{n\times(N+1)n}$ and $\mbQ_0,\mbQ_i\in \mathbb {R}^{(N+1)n\times(N+1)n}$.
We write
\begin{align}
&|\mbx_0-\Gamma_{0} \mbx^{(N)}-\eta_{0}|_{Q_{0}}^2 = \mbx^T \mbQ_{0} \mbx -2 \mbx^T \mbK_{0}^TQ_{0}\eta_{0} \nonumber\\
& \hskip 3.5cm + \eta_{0}^T Q_{0}\eta_{0},\label{sym}\\
&|\mbx_i-\Gamma_{1} \mbx_0-\Gamma_{2} \mbx^{(N)}-\eta|_{Q}^2 =
\mbx^T \mbQ_{i} \mbx - 2\mbx^T
\mbK_{i}^TQ\eta \nonumber\\
&\hskip 3.5cm+ \eta^T Q\eta, \quad 1\le i\le N.\label{sym2}
    \end{align}
 We may write $V_j(T,\mbx)$, $ 0\le j\le N$, in a similar form.

    We substitute \eqref{1st}  into \eqref{DE2} and derive the equation systems:
\begin{align}\label{DE3_P0}
\begin{cases}
\dot{\mbP}_0(t) =  - \big({\mbP}_0\widehat{\mbA}+\widehat{\mbA}^T \mathbb{\mbP}_0\big)+
                   {\mbP}_0{\mbB}_0 R_0^{-1} {\mbB}_0^T{\mbP}_0 \\
\qquad\quad        + \Big({\mbP}_0\sum_{k=1}^N {\mbB}_k R^{-1} {\mbB}_k^T {\mbP}_k\\
    \qquad\quad     +\sum_{k=1}^N {\mbP}_k{\mbB}_k R^{-1} {\mbB}_k^T {\mbP}_0\Big) -{\mbQ}_{0}   , \\
 {\mbP}_0(T) = {\mbQ}_{0f},
 \end{cases}
\end{align}
\begin{align}\label{DE3_S0}
\begin{cases}
\dot{{\mbS}}_0(t) = - \widehat{{\mbA}}^T {\mbS}_0  +  {\mbP}_0{\mbB}_0 R_0^{-1} {\mbB}_0^T {\mbS}_0\\
                            \qquad\quad    + \sum_{k=1}^N{\mbP}_k {\mbB}_k R^{-1} {\mbB}_k^T {\mbS}_0\\
                             \qquad\quad   + {\mbP}_0\sum_{k=1}^N {\mbB}_k R^{-1}{\mbB}_k^T {\mbS}_k
   +\mbK_0^T Q_{0}\eta_{0} , \\
{\mbS}_0(T)= -\mbK_{0f}^T Q_{0f}\eta_{0f},
\end{cases}
\end{align}
\begin{align}\label{DE3_gamma0}
\begin{cases}
\dot{\mbr}_0(t) = {\mbS}_0^T{\mbB}_0 R_0^{-1} {\mbB}_0^T{\mbS}_0
                    + 2{\mbS}_0^T\sum_{k=1}^N {\mbB}_k
                    R^{-1}{\mbB}_k^T {\mbS}_k \\
                   \qquad\quad  - \eta_{0}^T Q_{0} \eta_{0}
                    -\mbox{Tr}\big(\widehat{\mbD}^T
                          {\mbP}_0\widehat{\mbD}\big), \\
 \mbr_0(T)=\eta_{0f} ^T Q_{0f} \eta_{0f}.
\end{cases}
\end{align}

  By \eqref{DE3} and \eqref{1st}, we derive the equation systems:
\begin{align}\label{DE3_P}
\begin{cases}
\dot{\mbP}_i(t) =  - \big({\mbP}_i\widehat{\mbA}+\widehat{\mbA}^T \mathbb{\mbP}_i\big)- {\mbP}_i{\mbB}_i R^{-1} {\mbB}_i^T{\mbP}_i\\
    \qquad\quad    + \big({\mbP}_i{\mbB}_0 R_0^{-1} {\mbB}_0^T{\mbP}_0+
                           {\mbP}_0{\mbB}_0 R_0^{-1} {\mbB}_0^T{\mbP}_i\big)\\
    \qquad\quad   + \Big({\mbP}_i\sum_{k=1}^N {\mbB}_k R^{-1} {\mbB}_k^T {\mbP}_k+
 \\ \qquad\qquad
            \sum_{k=1}^N {\mbP}_k{\mbB}_k R^{-1} {\mbB}_k^T {\mbP}_i\Big)-{\mbQ}_{i}  \\
 {\mbP}_i(T) = {\mbQ}_{if},\qquad 1\leq i \leq N,
 \end{cases}
\end{align}
\begin{align}\label{DE3_S}
\begin{cases}
\dot{{\mbS}}_i(t) = - \widehat{{\mbA}}^T {\mbS}_i  +  {\mbP}_0{\mbB}_0 R_0^{-1} {{\mbB}_0}^T {\mbS}_i\\
          \qquad\quad    + {\mbP}_i{\mbB}_0 R_0^{-1} {\mbB}_0^T {\mbS}_0    
            -{\mbP}_i{\mbB}_i R^{-1} {\mbB}_i^T {\mbS}_i\\
          \qquad\quad    + \sum_{k=1}^N{\mbP}_k {\mbB}_k R^{-1} {\mbB}_k^T {\mbS}_i \\
          \qquad\quad
              + {\mbP}_i\sum_{k=1}^N {\mbB}_k R^{-1}{\mbB}_k^T {\mbS}_k   
              +\mbK_{i}^T Q \eta , \\
{\mbS}_i(T)= -\mbK_{if}^T Q_{f}\eta_{f}, \qquad 1\leq i \leq N,
\end{cases}
\end{align}
\begin{align}\label{DE3_gamma}
\begin{cases}
\dot{\mbr}_i(t) =  2{\mbS}_i^T{\mbB}_0 R_0^{-1} {\mbB}_0^T{\mbS}_0  \\
                   \qquad\quad  + 2{{\mbS}_i}^T\sum_{k=1}^N {\mbB}_k R^{-1}{\mbB}_k^T {\mbS}_k\\
                   \qquad\quad  - {\mbS}_i^T{\mbB}_i R^{-1} {\mbB}_i^T{\mbS}_i
                     - \eta^T Q \eta -\mbox{Tr}\big(\widehat{\mbD}^T
                          {\mbP}_i\widehat{\mbD}\big), \\
 \mbr_i(T)=\eta_{f} ^T Q_{f} \eta_{f}, \qquad 1\leq i \leq N.
\end{cases}
\end{align}

\begin{remark}\label{remark:P}
If \eqref{DE3_P0} and \eqref{DE3_P} have a solution $(\mbP_0, \cdots, \mbP_N)$  on $[\tau,T]\subseteq [0,T]$, such a solution is unique due to the local Lipschitz continuity of the vector field; see {\rm(Hale, 1969)}. The ODE guarantees  each $\mbP_j$, $0\le j\le N $, to be symmetric.
If \eqref{DE3_P0} and \eqref{DE3_P} have a unique solution $(\mbP_0,\cdots ,\mbP_N)$ on $[0, T]$, then we can uniquely solve $(\mbS_0, \cdots, \mbS_N)$ and $(\mbr_0, \cdots, \mbr_N)$.
\end{remark}


\begin{lemma}\label{theorem:Nash}
Suppose that \eqref{DE3_P0} and \eqref{DE3_P} have a unique solution $(\mbP_0,\cdots ,\mbP_N)$ on $[0,T]$. Then we can uniquely solve \eqref{DE3_S0}, \eqref{DE3_gamma0}, \eqref{DE3_S}, \eqref{DE3_gamma}, and the Nash game of $N+1$ players has a set of feedback Nash strategies given by
\begin{align}
&\hat u_0(t)=-R_0^{-1}\mbB_0^{T}(\mbP_0X(t)+\mbS_0),\label{Nashu0} \\
&\hat u_i(t)=-R^{-1} \mbB_i^T (\mbP_i X(t) +\mbS_i), \quad 1\le i\le N. \label{Nashui}
\end{align}
\end{lemma}
\proof This lemma follows the standard results in (Basar and Olsder, 1999, Theorem 6.16, Corollary 6.5).
\endproof

By Lemma \ref{theorem:Nash} and Remark \ref{remark:P},   the solution of the feedback Nash strategies completely reduces to the study  of \eqref{DE3_P0} and \eqref{DE3_P}.

\section{Asymptotic solvability}
\label{sec:asy}

Define the $(N+1)n\times (N+1)n$ identity matrix
\begin{align*}
I_{(N+1)n}=\begin{bmatrix}
I_n   &0    &\cdots        &0    \\
0  &I_n    &\cdots         &0     \\
\vdots  &\vdots     &\ddots        &\vdots     \\
0  &0    &\cdots         &I_n
\end{bmatrix}.
\end{align*}
For $1\le i\ne j\le N+1$,
exchanging the $i$th and $j$th rows of submatrices in $I_{(N+1)n}$,  let $J_{ij}$ denote the resulting matrix. For instance, we have
\begin{align*}
J_{12}=\begin{bmatrix}
0   &I_n    &\cdots        &0    \\
I_n  &0    &\cdots         &0     \\
\vdots  &\vdots     &\ddots        &\vdots     \\
0  &0    &\cdots         &I_n
\end{bmatrix}.
\end{align*}
It is easy to check that $J_{ij}^T=J_{ij}^{-1}=J_{ij}$.

\begin{theorem}\label{theorem:Prep}
$\mbP_0(t)$ and $\mbP_1(t)$ have the representation:
\begin{align}
&{\mbP}_0(t)=\begin{bmatrix}
\Pi_1^0 &\Pi_2^0 &\Pi_2^0&\cdots &\Pi_2^0 \\
{\Pi_2^0}^T &\Pi_3^0 &\Pi_3^0&\cdots &\Pi_3^0\\
{\Pi_2^0}^T &\Pi_3^0&\Pi_3^0 &\cdots &\Pi_3^0\\
\vdots      & \vdots      &  \vdots     &\ddots    &\vdots \\
{\Pi_2^0}^T &\Pi_3^0 &\Pi_3^0&\cdots &\Pi_3^0
\end{bmatrix}, \label{P0_matrix}
\end{align}
and
\begin{align}
&{\mbP}_1(t)=\begin{bmatrix}
\Pi_0   &\Pi_a   &\Pi_b  &\cdots   &\Pi_b \\
\Pi_a^T &\Pi_1   &\Pi_2  &\cdots   &\Pi_2\\
\Pi_b^T &\Pi_2^T &\Pi_3  &\cdots   &\Pi_3\\
\vdots     & \vdots    &  \vdots  &\ddots   &\vdots \\
\Pi_b^T &\Pi_2^T &\Pi_3  &\cdots   &\Pi_3
\end{bmatrix},\label{P1_matrix}
\end{align}
where each submatrix depends on $t$ and is $n\times n$.
Moreover, ${\mbP}_i(t)=J_{2,i+1}^T{\mbP}_1(t)J_{2,i+1}$ for $i\geq2$.
\end{theorem}

\proof See Appendix A.\hfill $\Box$

\begin{definition}\label{definition:as}
The sequence of  Nash games \eqref{stateX0}-\eqref{costJi} has asymptotic solvability if there exists $N_0$ such that for all $N\ge N_0$,
 $(\mbP_0,\cdots,\mbP_N)$ has a solution on $[0,T]$ and
 \begin{align}
 \sup_{N\ge N_0, 0\leq t\leq T} \left({\|\mbP_0(t)\|}_{l_1}+
{\|\mbP_1(t)\|}_{l_1}\right)<\infty. \label{P0P1}
\end{align}\hfill $\Box$
 \end{definition}

  Note that  \eqref{P0P1} is equivalent to
\begin{align}
\sup_{N\ge N_0, 0\leq t\leq T} &\left[|\Pi_1^0(t)|+N|\Pi_2^0(t)|+N^2|\Pi_3^0(t)|
\right]<\infty,\label{main_con0}\\
\sup_{N\ge N_0, 0\leq t\leq T} &\big[|\Pi_0(t)|+|\Pi_1(t)|+|\Pi_a(t)|+N|\Pi_b(t)|\nonumber\\
&+N|\Pi_2(t)|+N^2|\Pi_3(t)|\big]<\infty.\label{main_con}
\end{align}

   Denote $$M_0=B_0R_0^{-1}B_0^T,\quad M=BR^{-1}B^T.$$

Define
\begin{align}\label{L0P}
\left\{\begin{array}{lll}
\Lambda_1^{0N} = \Pi_1^0 ,&\Lambda_2^{0N} = N\Pi_2^0 ,\ &\Lambda_3^{0N} = N^2\Pi_3^0 ,\\
\Lambda_0^{N} = \Pi_0 ,&\Lambda_1^{N} = \Pi_1 ,&\Lambda_2^{N} = N\Pi_2 ,\\
\Lambda_3^{N} = N^2\Pi_3 ,\ &\Lambda_a^{N} = \Pi_a ,&\Lambda_b^{N} = N\Pi_b.
\end{array}\right.
\end{align}
For \eqref{DE3_P0} and \eqref{DE3_P}
we write the ODE system for the set of variables $(\Lambda_1^{0N}, \Lambda_2^{0N}, \dots, \Lambda_b^N)$ in \eqref{L0P}; see Appendix B. The following ODE system is obtained as the limit of the above ODE system with respect to $N$:
\begin{align}
\left\{
\begin{array}{l}
 \dot{\Lambda}_1^{0} =  \Lambda_1^{0}M_0 \Lambda_1^{0}-(\Lambda_1^{0}A_0+A_0^T\Lambda_1^{0}) \\
\quad\quad +\Lambda_2^{0}(M\Lambda_a^T-G)+(\Lambda_a M-G^T){\Lambda_2^{0}}^T-Q_0, \\
\dot{\Lambda}_2^{0} =( \Lambda_1^{0}M_0 -A_0^T) \Lambda_2^{0}+\Lambda_2^{0}(M(\Lambda_1+\Lambda_2)-A-F) \\
\quad\quad-\Lambda_1^{0}F_0+ (\Lambda_a M-G^T)\Lambda_3^{0}+Q_0\Gamma_0,\\
\dot{\Lambda}_3^{0} = {\Lambda_2^{0}}^TM_0 \Lambda_2^{0}-{\Lambda_2^{0}}^TF_0-F_0^T\Lambda_2^{0} \\
\quad\quad +\Lambda_3^{0}(M(\Lambda_1+\Lambda_2)-A-F) \\
\quad\quad +((\Lambda_1+\Lambda_2^T)M-A^T -F^T)\Lambda_3^{0}-\Gamma_0^TQ_0\Gamma_0,\\
 \dot{\Lambda}_0    = \Lambda_a M\Lambda_a^T-\Lambda_bG-G^T \Lambda_b^T \\
\quad\quad +\Lambda_0(M_0\Lambda_1^{0}-A_0)+(\Lambda_1^{0} M_0-A_0^T) \Lambda_0 \\
\quad\quad -\Lambda_a(G-M \Lambda_b^T)-(G^T-\Lambda_bM)
\Lambda_a^T-\Gamma_1^TQ\Gamma_1, \\
\dot{\Lambda}_1 = \Lambda_1M\Lambda_1-\Lambda_1A-A^T\Lambda_1-Q,\\
 \dot{\Lambda}_2 = \Lambda_a^T(M_0\Lambda_2^{0}-F_0)-\Lambda_1F+(\Lambda_1M-A^T)\Lambda_2 \\
\quad\quad+\Lambda_2(M(\Lambda_1+\Lambda_2)-A-F)+Q\Gamma_2,\\
\dot{\Lambda}_3 = \Lambda_b^TM_0\Lambda_2^{0}+{\Lambda_2^{0}}^TM_0\Lambda_b+\Lambda_2^TM\Lambda_2 \\
\quad\quad -\Lambda_b^TF_0-F_0^T\Lambda_b-\Lambda_2^TF-F^T\Lambda_2  \\
\quad\quad +\Lambda_3(M(\Lambda_1+\Lambda_2)-A-F) \\
\quad\quad +((\Lambda_1+\Lambda_2^T) M-A^T-F^T)\Lambda_3- \Gamma_2^T Q\Gamma_2,\\
\dot{\Lambda}_a = (\Lambda_1^{0}M_0-A_0^T)\Lambda_a+\Lambda_a(M\Lambda_1-A) \\
\quad\quad-G^T\Lambda_1+(\Lambda_aM-G^T)\Lambda_2^T+\Gamma_1^T Q,\\
 \dot{\Lambda}_b = \Lambda_0M_0\Lambda_2^{0}+(\Lambda_aM-G^T)(\Lambda_2+\Lambda_3)-\Lambda_0 F_0 \\
\quad\quad
-\Lambda_a F  +\Lambda_b(M(\Lambda_1+\Lambda_2)-A-F)\\
\quad\quad  +(\Lambda_1^{0}M_0-A_0^T)\Lambda_b -\Gamma_1^T Q\Gamma_2,
\end{array}
\right.\label{eqn1*}
\end{align}
where the terminal conditions are
\begin{align}
\begin{cases}
 \Lambda_1^0(T)=Q_{0f},\quad
 \Lambda_2^0(T)=-Q_{0f}\Gamma_{0f},\nonumber\\
 \Lambda_3^0(T)=\Gamma_{0f}^TQ_{0f}\Gamma_{0f},\nonumber  \\
\Lambda_0 (T)=\Gamma_{1f}^TQ_f\Gamma_{1f},\quad
\Lambda_1 (T)=Q_{f},\nonumber \\
  \Lambda_2 (T)=- Q_f\Gamma_{2f},\quad
   \Lambda_3 (T)=\Gamma_{2f}^T Q_f\Gamma_{2f},\nonumber\\
 \Lambda_a (T)=-\Gamma_{1f}^T Q_f ,\quad
\Lambda_b (T)=\Gamma_{1f}^T Q_f\Gamma_{2f}.\nonumber  \\
\end{cases}
\end{align}

\begin{theorem}\label{theorem:AS}
The sequence of games in \eqref{stateX0}-\eqref{costJi} has asymptotic solvability if and only if \eqref{eqn1*} has a  solution on $[0,T]$.
\end{theorem}

\proof See Appendix B.\hfill$\Box$

Due to the quadratic terms in its right hand sides, we call
\eqref{eqn1*} a system of Riccati ODEs.
As it turns out later in Section \ref{sec:two}, this set of solution functions  can be interpreted according to two optimal  control problems.

\section{Equilibrium costs and decentralized control}
\label{sec:cos}

For this section, we assume \eqref{eqn1*} has  a solution on $[0,T]$.
Therefore there exists $N_0>0$ such that for all $N\ge N_0$,
\eqref{DE3_P0} and \eqref{DE3_P} have a solution  $(\mbP_0, \cdots, \mbP_N) $ on $[0,T]$.
\begin{proposition}\label{prop:S}
Let $(\mbS_{0},\cdots,\mbS_{N})$ be the solution of \eqref{DE3_S0} and \eqref{DE3_S}. We have the representation
\begin{align} 
&\mbS_{0}(t)= [{\theta_0^0}^T,{\theta_1^0}^T,\cdots,{\theta_1^0}^T]^T, \nonumber\\
&\mbS_{i}(t)=[{\theta^T_0},{\theta^T_2},\cdots,{\theta^T_2},{\theta^T_1},{\theta^T_2},\cdots,{\theta^T_2}]^T,\quad 1\leq i\leq N,\nonumber
\end{align}
where each vector of ${\theta_0^0(t),~\theta_1^0(t),~\theta_0(t),~\theta_1(t),~\theta_2(t)}$ is in $\mathbb{R}^n$ and
$\theta_1^T$ is the $(i+1)$th component of $\mbS_{i}$.
\end{proposition}

\proof  The method is similar to proving Theorem \ref{theorem:Prep}, and we omit the detail. \endproof

Define
\begin{align}
&\alpha_0^{0N} = \theta_0^0 ,\quad \alpha_1^{0N} = N\theta_1^0 ,\nonumber\\
&\alpha_0^{N} = \theta_0 ,\quad \alpha_1^{N} = \theta_1 ,\quad
\alpha_2^{N} = N\theta_2. \nonumber
\end{align}
We derive a set of ODEs for $(\alpha_0^{0N},\alpha_1^{0N},\alpha_0^{N},\alpha_1^{N},\alpha_2^{N})$; see Appendix C.
By taking the limit form of these equations with respect to $N$, we introduce the ODE system:
\begin{align}
\left\{\begin{array}{l}
\dot{\alpha}_0^{0} = (\Lambda_1^{0}
M_0-A_0^T)\alpha_0^{0}+(\Lambda_aM-G^T)\alpha_1^{0} \nonumber\\
\quad\quad +\Lambda_2^{0}M\alpha_1+Q_0\eta_0, \\
\dot{\alpha}_1^{0} =({\Lambda_2^{0}}^TM_0 -F_0^T)\alpha_0^{0}+
(\Lambda_1 M+ \Lambda_2 ^TM-A^T-F^T)\alpha_1^{0}
 \nonumber\\
\quad\quad
+\Lambda_3^{0}M\alpha_1-\Gamma_0^TQ_0\eta_0, \\
\dot{\alpha}_0 = (\Lambda_1^{0}M_0-A_0^T)\alpha_0+((\Lambda_a+\Lambda_b)M-G^T)\alpha_1
 \nonumber\\
\quad\quad +\Lambda_0M_0\alpha_0^{0}+(\Lambda_aM-G^T)\alpha_2-\Gamma_1^TQ\eta, \\
\dot{\alpha}_1 = ((\Lambda_1+\Lambda_2)M-A^T)\alpha_1+\Lambda_a^TM_0\alpha_0^0+Q\eta,  \\
\dot{\alpha}_2 = ({\Lambda_2^0}^TM_0-F_0^T)\alpha_0+ ((\Lambda_1+\Lambda_2 ^T) M -A^T-F^T)\alpha_2\nonumber\\
\quad\quad+\Lambda_b^TM_0 \alpha_0^0  +((\Lambda_2 ^T+\Lambda_3)M-F^T)\alpha_1-\Gamma_2^TQ\eta,
\end{array}\right.
\end{align}
where
\begin{align}
\begin{cases}
 \alpha _0^{0 }(T)=-Q_{0f}\eta_{0f}, \quad
 \alpha _1^{0 }(T)=\Gamma_{0f}^TQ_{0f}\eta_{0f},\nonumber\\
  \alpha _0 (T)= \Gamma_{1f}^TQ_f\eta_f, \quad
\alpha_1 (T)= - Q_f\eta_f, \nonumber\\
  \alpha _2 (T)=\Gamma_{2f}^TQ_f\eta_f. \nonumber
\end{cases}
\end{align}
After \eqref{eqn1*} is solved, $(\alpha_0^0,\cdots,\alpha_2)$ satisfies a linear ODE system and can be uniquely solved on $[0,T]$.

\begin{proposition}\label{prop:thal}
We have
\begin{align}
&\sup_{0\le t\le T}\{|\theta_0^{0 }(t)-\alpha_0^{0}(t)|+|N\theta_1^{0}(t)-\alpha_1^{0}(t)|
+ |\theta_0 (t)-\alpha_0 (t)|\nonumber\\
&\qquad\qquad+ |\theta_1(t)-\alpha_1(t)|+|N\theta_2(t)-\alpha_2(t)|\}=O(\tfrac{1}{N}).\nonumber
\end{align}
\end{proposition}
\proof We consider the first ODE system for $(\Lambda_1^{0N},\cdots, \Lambda_b^{N} )$ and $(\alpha_0^{0N}, \cdots, \alpha_2^N)$, and the second ODE system for $(\Lambda_1^{0},\ldots, \Lambda_{b} )$ and $(\alpha_0^{0}, \cdots, \alpha_2)$. By (Huang and Zhou, 2018b, Theorem 4), we obtain the error bound.
  \endproof

In view of \eqref{DE3_gamma0}  and \eqref{DE3_gamma}, we obtain
\begin{align}
&\dot{\mbr}_0={\theta_0^0}^TM_0\theta_0^0+2N{\theta_1^0}^TM
\theta_1-\eta_0^TQ_0\eta_0\nonumber\\
&\quad\quad - {\rm Tr}(D_0^T\Pi_1^0D_0)-{\rm Tr}(ND^T\Pi_3^0D),\nonumber\\
&\dot{\mbr}_i =2\theta_0^TM_0\theta_0^0+2(N-1)\theta_2^TM\theta_1+\theta_1^TM \theta_1-\eta^TQ\eta\nonumber\\
&\quad\quad -{\rm Tr}(D_0^T\Pi_0D_0+D^T\Pi_1D)-(N-1){\rm Tr}(D^T\Pi_3D),\nonumber
\end{align}
where ${\mbr}_0(T)=\eta_{0f}^TQ_{0f}\eta_{0f} $
and ${\mbr}_i (T)=\eta_f^TQ_f\eta_f$. For $N\geq N_0$, we can uniquely solve ${\mbr}_0$ and ${\mbr}_i$. It is clear that ${\mbr}_i$ does not depend on $i$.
We rewrite
\begin{align}
&\dot{\mbr}_0=(\alpha_0^{0N})^TM_0\alpha_0^{0N}+2(\alpha_1^{0N})^TM\alpha_1^N
-\eta_0^TQ_0\eta_0 \nonumber \\
&\quad\quad-{\rm Tr}(D_0^T\Lambda_1^{0N}D_0)
-\tfrac{1}{N}{\rm Tr}(D^T\Lambda_3^{0N}D),\label{r0}\\
&\dot{\mbr}_i=2{\alpha_0^N}^TM_0\alpha_0^{0N}+2 {\alpha_2^N}^TM\alpha_1^N+{\alpha_1^N}^TM\alpha_1^N \nonumber\\
&\quad\quad  -\eta^TQ\eta-{\rm Tr}(D_0^T\Lambda_0^ND_0)-{\rm Tr}(D^T\Lambda_1^ND) \nonumber\\
&\quad\quad -\tfrac{2}{N}{\alpha_2^N}^TM\alpha_1^N-
(\tfrac{1}{N}-\tfrac{1}{N^2}){\rm Tr}(D^T\Lambda_3^ND).\label{ri}
\end{align} 
As the approximation of \eqref{r0} and \eqref{ri},  we  introduce the ODE system
\begin{align}
\left\{
\begin{array}{l}
\dot{\chi}_0={\alpha_0^{0}}^TM_0\alpha_0^{0}+2{\alpha_1^0}^TM\alpha_1-\eta_0^T
Q_0\eta_0-{\rm Tr}(D_0^T\Lambda_1^{0}D_0),\\
\dot{\chi} =2\alpha_0^TM_0\alpha_0^0 +2 \alpha_2 ^TM\alpha_1 +\alpha_1^TM\alpha_1 \\
\quad\quad -\eta^TQ\eta-{\rm Tr}(D_0^T\Lambda_0 D_0)-{\rm Tr}(D^T\Lambda_1 D),
\end{array}\right.\nonumber
\end{align}
where
$\chi_0(T)=\eta_{0f}^TQ_{0f}\eta_{0f}$ and $ 
\chi (T)=\eta_f ^TQ_f \eta_f$, and we solve $(\chi_0,~\chi)$ on $[0,~T]$.

\begin{proposition}
We have
\begin{align}
\sup_{0\le t\le T}\{|{\mbr}_0(t)-\chi_0 (t)|+|{\mbr}_i(t)-\chi (t)|\}=O(\tfrac{1}{N}).\nonumber
\end{align}
\end{proposition}
\proof  The proof is similar to that of Proposition \ref{prop:thal}.  \endproof

Assumption ({H}):
The initial states $X_1(0),X_2(0),\cdots$, are i.i.d. and $X_1(0)$ has mean $\mu$, and covariance $\Sigma$. In addition,   $X_0(0)$ has mean $\mu_0$ and covariance $\Sigma_0.$

 Denote the set of Nash  strategies
 $\hat{u}=(\hat{u}_0,\hat{u}_1,\cdots,\hat{u}_N)$ given by \eqref{Nashu0}-\eqref{Nashui}.
\begin{proposition} Under Assumption (H),
the costs under the set of strategies $\hat{u}$ have the asymptotic form
\begin{align}
&\lim_{N\to\infty}J_0(\hat{u})=\mu_0^T\Lambda_1^0(0)\mu_0
+2\mu_0^T\Lambda_2^0(0)\mu+\mu^T\Lambda_3^0\mu \nonumber\\
&\quad\quad +2(\mu_0^T{\alpha_0^0}(0)+\mu^T{\alpha_1^0}(0))+{\rm Tr}(\Lambda_1^0(0)\Sigma_0)+\chi_0(0), \nonumber
\end{align}
\begin{align}
&\lim_{N\to\infty}J_1(\hat{u})=\mu_0^T\Lambda_0(0)\mu_0
+2\mu^T\Lambda_a(0)\mu_0 +2\mu_0^T\Lambda_b(0)\mu\nonumber\\
&\quad\quad +\mu^T(\Lambda_1(0)
+2\Lambda_2(0) +\Lambda_3(0))\mu \nonumber \\
&\quad\quad +2[\mu_0^T\alpha_0(0) +\mu^T \alpha_1(0) +\mu^T\alpha_2(0) ]  \nonumber\\
&\quad\quad +{\rm Tr}(\Lambda_0(0)\Sigma_0+\Lambda_1(0)\Sigma)+\chi_1(0). \nonumber
\end{align}
\end{proposition}

\proof
Note that
\begin{align}
J_0(\hat{u})= E[{X^T(0)}\mbP_0(0)X(0)+2{{\mbS}^T_0(0)}X(0)]+\mbr_0(0), \nonumber
\end{align}
and
\begin{align}
&{X^T(0)}{\mbP}_0(0)X(0)={X^T_0(0)}\Pi_1^0(0)X_0(0)\nonumber\\
&\quad\quad  +2{X^T_0(0)}[\Pi_2^0(0),\cdots,\Pi_2^0(0)]X_{-0}(0) \nonumber\\
&\quad\quad +{X^T_{-0}(0)}
\begin{bmatrix}
\Pi_3^0(0) & \cdots & \Pi_3^0(0) \\
\vdots      & \ddots      &  \vdots \\
\Pi_3^0(0) & \cdots & \Pi_3^0(0) \\
\end{bmatrix}
 X_{-0}(0), \nonumber
\end{align}
where
$X_{-0}(0)=[X_1^T(0),
X_2^T(0),
\cdots,
X_N^T(0)]^T$.
Similarly we have
\begin{align}
J_1(\hat{u})=E[X^T(0)\mbP_1(0)X(0)+2\mbS^T_1(0)X(0)]+\mbr_1(0).\nonumber
\end{align}
We complete the proof by elementary computations and taking limits. \endproof

Substituting $\hat u_0$ and $\hat u_i$ into \eqref{stateX0} and \eqref{stateXi}, we have
\begin{align}
&dX_0=\big(A_0X_0-M_0(\Pi_1^0X_0+N\Pi_2^0X^{(N)}+\theta_0^0)\nonumber\\
&\quad\quad +F_0X^{(N)}\big)dt+D_0dW_0,\nonumber\\
& dX^{(N)}=\big(AX^{(N)}-M(\Pi_a^TX_0+\Pi_1X^{(N)}+(N-1)\Pi_2X^{(N)}\nonumber\\
&\qquad\qquad +\theta_1)+FX^{(N)}+GX_0\big)dt+\frac{1}{N}
\sum_{i=1}^{N}DdW_i. \label{xiN}
 \end{align}
When $N\to \infty$, we   obtain a limit form of the strategies
 \begin{align}
&\check u_0=-R_0^{-1}B_0^T (\Lambda_1^0X_0+\Lambda_2^0\overline
 X+\alpha_0^0), \label{u0}\\
& \check u_i=-R^{-1}B^T (\Lambda_a^T X_0+\Lambda_1X_i+\Lambda_2\overline X+\alpha_1) ,\label{ui}
\end{align}
where $\overline X$ is the infinite population limit of the  state average ${X}^{(N)}$ of the minor players in \eqref{xiN}.
For the $N+1$ player game, we replace $\overline X$ in
 the  strategies \eqref{u0}-\eqref{ui} by $\overline X^\dag$ and write the
closed-loop system of equations:
\begin{align}
dX_0(t)=&\big[A_0X_0-M_0(\Lambda_1^0X_0+\Lambda_2^0\overline X^\dag+\alpha_0^0)\nonumber\\
&\quad +F_0X^{(N)}\big]dt+D_0dW_0,\label{XX0}\\
dX_i(t)=  &\big[AX_i-M (\Lambda_a^T X_0+\Lambda_1X_i+\Lambda_2\overline X^\dag+\alpha_1)\nonumber \\
&+FX^{(N)}+GX_0\big]dt +DdW_i, \quad 1\le i\le N,  \label{XXI}  \\
d\overline {X}^\dag(t)=&[(A-M(\Lambda_1+\Lambda_2)+F )\overline {X}^\dag
 \nonumber\\
 &+(G-M \Lambda_a^T ) X_0-M\alpha_1]dt, \qquad t\ge 0,\nonumber
\end{align}
where $\overline X^\dag$ is generated by the $N+1$ players instead of an infinite population. Denote the strategies in \eqref{XX0}-\eqref{XXI} by $(\check u_0^\dag, \cdots, \check u_N^\dag)$.
Following the standard mean square error estimate of $|X^{(N)}-\overline X|$ for \eqref{XX0}-\eqref{XXI}, under assumption (H) we can show that $(\check u_0^\dag, \cdots, \check u_N^\dag)$ is an $\epsilon$-Nash equilibrium for the $N+1$ player game, where $\epsilon =O(1/\sqrt{N})$ and each player may use centralized state information $X(t)$; see related methods in (Huang, 2010).

By the re-scaling technique, we derive the mean field limits of the costs and strategies. The feasibility condition is determined by \eqref{eqn1*} directly based on the model parameters in \eqref{stateX0}-\eqref{costJi}. This is different from (Huang, 2010), where the existence condition is described in an augmented state space in $3n$ dimensions and imposes consistency requirements on $3n\times 3n$ matrices.

\section{The limiting control problems and best responses}
\label{sec:two}

For this section, we assume \eqref{eqn1*} has  a solution on $[0,T]$.

An interesting question is whether the above two limit strategies in \eqref{u0}-\eqref{ui} have the interpretation as best responses in appropriately constructed optimal control problems. Finding the best response of a single agent in an infinite population model has been a key step in the fixed point approach in mean field games; see (Huang, Caines, and Malham\'e, 2007; Huang, 2010). We introduce two optimal control problems.

Problem (P0): The dynamics are given by
\begin{align}
dX_0(t)=&\big(A_0X_0+B_0u_0+F_0\overline{X}\big)dt+DdW_0, \label{x0xb} \\
d\overline {X}(t)=&[(A-M(\Lambda_1+\Lambda_2)+F )\overline {X}
 \nonumber\\
&+(G-M \Lambda_a^T) X_0-M\alpha_1]dt,\label{state0Xi}
\end{align}
where $X_0(0)$ and $\overline{X}(0)=\mu_0$ are given.
Equation \eqref{state0Xi} may be viewed as the limit of \eqref{xiN}
but now $X_0$ is indirectly controlled by $u_0$ in \eqref{x0xb}.
The cost  is
\begin{align}
\overline{J}_0(u_0) =\  & E\int_0^T \Big(  |X_0(t) -\Gamma_0 \overline{X}(t) -\eta_0|_{Q_0}^2+|u_0(t)|^2_{R_0}\Big)dt \nonumber\\
&+E|X_0(T)-\Gamma_{0f} \overline{X}(T)-\eta_{0f}|_{Q_{0f}}^2.\nonumber
\end{align}

Problem (P1): The dynamics are given by
\begin{align}
dX_1(t)=&\big(A X_1+B u_1+F \overline{X}+GX_0\big)dt+DdW_1,\nonumber \\
d {X_0}(t)=&[A_0X_0-M_0(\Lambda_1^0X_0+\Lambda_2^0\overline X+\alpha_0^0)]dt \nonumber\\
&+F_0\overline X dt+D_0dW_0 , \label{x0xbar}\\
d\overline {X}(t)=&[(A-M(\Lambda_1+\Lambda_2)+F )\overline {X}
  \nonumber\\
 &+(G-M \Lambda_a^T) X_0-M\alpha_1]dt, \nonumber
\end{align}
where $X_1(0),~X_0(0)$ and $\overline X(0)=\mu_0$ are given. The  notation  $(X_0,\overline X)$ is reused  in Problem (P1), where $u_0$ has  taken a specific form in \eqref{x0xbar}. Equation \eqref{x0xbar}  can be viewed as a limit form of \eqref{XX0} when $N\to \infty$.  Since the two problems will be solved separately, this should cause no risk of confusion.
The cost is
\begin{align}
{\overline J}_1(u_1) =\  & E\int_0^T \Big(  |X_1(t)-\Gamma_1 X_0(t)-\Gamma_2 \overline{X}(t)-\eta|_{Q}^2+|u_1|^2_{R}\Big)dt \nonumber\\
&+E|X_1(T)-\Gamma_{1f} X_0(T)-\Gamma_{2f}
\overline X (T)-\eta_{f}|_{Q_{f}}^2, \nonumber
\end{align}
 Since $R_0>0, ~ Q_0,~Q_{0,f}\geq 0, ~R>0, ~ Q,~Q_f\geq0$, both Problem (P0) and Problem (P1) can be solved. The resulting optimal control laws will also be called best responses.

Below, we start with the solution of Problem (P0). Denote
\begin{align}
&\mathbb{A}_0=\left[
\begin{array}{cc}
{A}_{0} & F_0 \\
  G-M\Lambda_a^T & A+F-M(\Lambda_1+\Lambda_2)
  \end{array}
  \right],
  \quad
 \mathbb{B}_0=\left[\begin{matrix}
  B_0\\
    0\end{matrix}\right],
  \nonumber\\
&  \mathbb{Q}_0=\left[
\begin{array}{c}
I_n \\
-\Gamma_0^T
\end{array}\right]
{Q}_0 [I_n,-\Gamma_0].\nonumber
\end{align}
Then we have
\begin{align}
&d\left[
\begin{array}{c}
X_0 \\
\overline{X}
\end{array}\right]
=\mathbb{A}_0
\left[\begin{array}{c}
 X_0  \\
 \overline{X}
\end{array}\right]dt+\mathbb B_0 u_0dt-
\left[\begin{array}{c}
 0  \\
 M\alpha_1
\end{array}\right]dt +
\left[\begin{array}{c}
 D_0dW_0  \\
 0
\end{array}\right].\nonumber
\end{align}

We further write
\begin{align}
&|X_0-\Gamma_0 \overline{X}-\eta_0|_{Q_0}^2\nonumber\\
=\ &(X_0^T,\overline{X}^T)
{\mathbb Q}_0
\left[\begin{array}{c}
 X_0  \\
 \overline{X}
\end{array}\right]
-2\eta_0^TQ_0[I_n,-\Gamma_0]\left[\begin{array}{c}
 X_0  \\
 \overline{X}
\end{array}\right]\nonumber\\
&+\eta_0^TQ_0\eta_0.\nonumber
\end{align}
By dynamic programming for the optimal control problem (P0), we introduce
\begin{align}
&\mathbb{\dot P}_0=-\mathbb A_0^T\mathbb P_0-\mathbb P_0\mathbb A_0+\mathbb P_0\mathbb B_0R_0^{-1}\mathbb B_0^T\mathbb P_0-\mathbb Q_0,\label{P0}\\
&\mathbb{\dot S}_0=-\mathbb A_0^T\mathbb S_0 + \mathbb P_0 \mathbb B_0R_0^{-1}\mathbb B_0^T\mathbb S_0  \nonumber\\
&\quad\quad +\mathbb P_0
\left[\begin{array}{c}
0\\
M\alpha_1\end{array}\right]+
\left[\begin{array}{c}
I_n \\
-\Gamma_0^T
\end{array}\right]Q_0\eta_0,\label{S0}
\end{align}
where the terminal condition can be determined as
\begin{align}
&\mathbb{P}_0(T)=\left[\begin{array}{c}
I_n \\
-\Gamma_{0f}^T
\end{array}\right]Q_{0f}[I_n,-\Gamma_{0f}],
\nonumber\\
&{\mathbb S}_0(T) = -
\begin{bmatrix}
I_n \\
-\Gamma_{0f}^T
\end{bmatrix}Q_{0f}\eta_{0f}.\nonumber
\end{align}
We uniquely solve ${\mathbb P}_0$ and ${\mathbb S}_0$ on $[0,T]$.
Note that $\mathbb P_0$ is a $2n\times 2n$ matrix.
The optimal control law is
$$
u_0^*(t)=-R_0^{-1}{\mathbb B}_0^T({\mathbb P}_0 [X_0^T(t), {\overline X}^T(t)]^T+{\mathbb S}_0).
$$
Denote
\begin{align}
&\mathbb P_0=\left[\begin{array}{ll}
\Phi_{1}^0&\Phi_{2}^0\\
{\Phi_2^0}^T&\Phi_{3}^0
\end{array}\right],~~\Phi_{k}^0(t)\in\mathbb{R}^{ n\times n},\quad
\mathbb S_0=\left[\begin{array}{lll}
\beta_0^0\\
\beta_1^0
\end{array}\right],\nonumber
\end{align}
where $\Phi_{1}^0$ and $\Phi_{3}^0$ are symmetric.
Then by \eqref{P0}, we derive the ODE system:
\begin{align}
\left\{
\begin{array}{l}
\dot{\Phi}_{1}^0=\Phi_{1}^0M_0\Phi_{1}^0-A_0^T\Phi_{1}^0
-\Phi_{1}^0A_0-(G^T-\Lambda_aM){\Phi_{2}^0}^T\nonumber\\
\quad\quad~~ -\Phi_{2}^0(G-M\Lambda_a^T)-Q_0,\\
\dot{\Phi}_{2}^0=-A_0^T\Phi_{2}^0-(G^T-\Lambda_a M)\Phi_{3}^0-\Phi_{1}^0F_0\nonumber\\
\quad\quad~~ -\Phi_{2}^0(A+F-M(\Lambda_1+\Lambda_2))\nonumber\\
\quad\quad~~ +{\Phi_{1}^0} M_0 \Phi_{2}^0+Q_0\Gamma_0,\\
\dot{\Phi}_{3}^0={\Phi_{2}^0}^T M_0 \Phi_{2}^0-F_0^T\Phi_{2}^0-{\Phi_{2}^0}^TF_0\nonumber\\
\qquad\quad-(A^T+F^T-(\Lambda_1+\Lambda_2^T)M)\Phi_{3}^0\nonumber\\
\quad\quad~~ -\Phi_{3}^0(A+F-M(\Lambda_1+\Lambda_2)) -\Gamma_0^TQ_0\Gamma_0,
\end{array}\right.
\end{align}
where
\begin{align}
\Phi_{1}^0(T)=Q_{0f}, \quad \Phi_{2}^0(T)=-Q_{0f}\Gamma_{0f},
\quad \Phi_{3}^0(T)=\Gamma_{0f}^TQ_{0f}\Gamma_{0f} .\nonumber
\end{align}
And by \eqref{S0},
\begin{align}
\left\{\begin{array}{l}
\dot{\beta}_0^0=(\Phi_{1}^0M_0-A_0^T)\beta_0^0-(G^T-\Lambda_a M)\beta_1^0\nonumber\\
\quad\quad~~ +\Phi_{2}^0M\alpha_1+Q_0\eta_0,\\
\dot{ \beta}_1^0=({\Phi_2^0}^TM_0 -F_0^T)\beta_0^0-(A^T+F^T-(\Lambda_1+\Lambda_2^T)M)\beta_1^0\nonumber\\
\quad\quad~~ +\Phi_{3}^0M\alpha_1-\Gamma_0^TQ_0\eta_0,
\end{array}\right.
\end{align}
where
\begin{align}
\beta_0^0(T)=-Q_{0f}\eta_{0f},\quad
\beta_1^0(T)=\Gamma_{0f}^TQ_{0f}\eta_{0f}.\nonumber
\end{align}
Finally, we rewrite $u_0^*$ as
$$
u_0^*(t)= -R_0^{-1}B_0^T(\Phi_1^0 X_0+\Phi_2^0\overline X+\beta_0^0).
$$

Now we give the solution of Problem (P1). Denote
\begin{align}
&\mathbb A=\left[
\begin{matrix}
&A_0-M_0\Lambda_1^0 &0 &F_0-M_0\Lambda_2^0\\
&G &A &F\\
 &G-M\Lambda_a^T &0 &A+F-M(\Lambda_1+\Lambda_2)
\end{matrix}
\right],\nonumber\\
&\mathbb B=\left[\begin{array}{ll}
0\\
B\\
0\end{array}\right],\qquad
f=-\left[
\begin{matrix}
M_0\alpha_0^0\\
0\\
M\alpha_1
\end{matrix}
\right],\nonumber\\
&\mathbb Q=[-\Gamma_1,I_n,-\Gamma_2]^TQ[-\Gamma_1,I_n,-\Gamma_2],\nonumber\\
&\mathbb Q_f=[-\Gamma_{1f},I_n,-\Gamma_{2f}]^TQ_f[-\Gamma_{1f},I_n,-\Gamma_{2f}].\nonumber
\end{align}

The dynamics can be given as
\begin{align}
\left[\begin{matrix}
dX_0\\
dX_1\\
d\overline X
\end{matrix}\right]=\mathbb A\left[
\begin{matrix}
X_0\\
X_1\\
\overline X
\end{matrix}
\right]dt+(\mathbb B u_1+f)dt+\left[\begin{matrix}
D_0dW_0\\
D_1dW_1\\
0
\end{matrix}\right]. \nonumber
\end{align}

By dynamic programming,
we introduce the two ODEs:
\begin{align}
&\mathbb{\dot P}=-\mathbb A^T\mathbb P-\mathbb P\mathbb A+\mathbb P\mathbb B R ^{-1}\mathbb B ^T\mathbb P -\mathbb Q, \label{P}\\
&\mathbb{\dot S}=-\mathbb A ^T\mathbb S  + \mathbb P  \mathbb B R ^{-1}\mathbb B ^T\mathbb S \nonumber\\
&\quad\quad +\mathbb P
\left[\begin{array}{c}
M_0\alpha_0^0\\
0\\
M\alpha_1\end{array}\right]+
\left[\begin{array}{c}
-\Gamma_1^T\\
I_n \\
-\Gamma_2^T
\end{array}\right]Q\eta,\label{S}
\end{align}
where
$$
{\mathbb P}(T) = {\mathbb Q}_f,
\quad {\mathbb S}(T)=-
\begin{bmatrix}
-\Gamma_{1f}^T\\
I_n\\
-\Gamma_{2f}^T
\end{bmatrix} Q_f \eta_f.
$$
We uniquely solve ${\mathbb P}$ and ${\mathbb S}$ on $[0,T]$.
Denote
\begin{align}
&\mathbb P =\left[\begin{array}{ccc}
\Phi_{0},&\Phi_{a},&\Phi_{b}\\
{\Phi_{a}^T},&\Phi_{1},&\Phi_{2}\\
{\Phi_{b}^T},&\Phi_{2}^T,&\Phi_{3}
\end{array}\right],\quad \mathbb S=\left[\begin{array}{cc}
\beta_0\\
\beta_1\\
\beta_2
\end{array}\right],\nonumber
\end{align}
where $\Phi_{0}$, $\Phi_1$ and $\Phi_3$ are symmetric.
By \eqref{P}, we derive the ODE system:
\begin{align}
\left\{
\begin{array}{l}
 \dot{\Phi}_{0}=\Phi_{a}M\Phi_{a}^T-G^T\Phi_{a}^T-\Phi_{a}G-\Gamma_1^TQ\Gamma_1\nonumber\\
\quad\quad -(A_0^T-{\Lambda_1^0} M_0)\Phi_{0}-\Phi_{0}(A_0- M_0{\Lambda_1^0}) \nonumber\\
\quad\quad -(G^T-\Lambda_aM)\Phi_{b}^T-\Phi_{b}(G-M\Lambda_a^T),  \\
 \dot{\Phi}_{1}=\Phi_{1}M\Phi_{1}-A^T\Phi_{1}-\Phi_{1}A-Q , \\
 \dot{\Phi}_{2}=-A^T\Phi_{2}-\Phi_{1}F-\Phi_{a}^T(F_0-M_0\Lambda_2^0) \nonumber\\
\quad\quad -\Phi_{2}(A+F-M(\Lambda_1+\Lambda_2))+\Phi_{1}M\Phi_{2}\nonumber\\
\quad\quad +Q\Gamma_2,  \\
 \dot{\Phi}_{3}=\Phi_{2}^TM\Phi_{2}-F^T\Phi_{2}-\Phi_{2}^TF \nonumber\\
\quad\quad  -(F_0^T-{\Lambda_2^0}^T M_0)\Phi_{b}-\Phi_{b}^T(F_0- M_0{\Lambda_2^0})\nonumber\\
\quad\quad  -(A^T+F^T-(\Lambda_1+\Lambda_2^T)M)\Phi_{3} \nonumber\\
\quad\quad-\Phi_{3}(A+F -M(\Lambda_1+\Lambda_2))-\Gamma_2^TQ\Gamma_2, \\
 \dot{\Phi}_{a}=-\Phi_{a}A-G^T\Phi_{1}-(A_0^T-\Lambda_1^0M_0)\Phi_a\nonumber\\
\quad\quad -(G^T-\Lambda_a M)\Phi_2^T+\Phi_{a}M\Phi_{1}+\Gamma_1^T Q,  \\
 \dot{\Phi}_{b}=-\Phi_{0}(F_0-M_0\Lambda_2^0)-G^T\Phi_{2}-\Phi_{a}F+\Phi_{a}M\Phi_{2} \nonumber\\
\quad\quad -(A_0^T-{\Lambda_1^0} M_0)\Phi_{b}-(G^T-\Lambda_aM)\Phi_{3}  \nonumber\\
\quad\quad  -\Phi_{b}(A+F-M  (\Lambda_1+\Lambda_2))-\Gamma_1^TQ\Gamma_2,
\end{array}\right.
\end{align}
where \begin{align}
& \Phi_{0}(T)=\Gamma_{1f}^TQ_f\Gamma_{1f}, \quad
 \Phi_{1}(T)=Q_f, \quad
 \Phi_{2}(T)=-Q_f\Gamma_{2f}, \nonumber\\
& \Phi_{3}(T)=\Gamma_{2f}^TQ_f\Gamma_{2f},  \quad
 \Phi_{a}(T)=-\Gamma_{1f}^T Q_f, \nonumber\\
& \Phi_{b}(T)=\Gamma_{1f}^TQ_f\Gamma_{2f}.  \nonumber
\end{align}
Now by \eqref{S}, we derive
\begin{align}
\left\{
\begin{array}{l}
 \dot{\beta}_0=(\Phi_{a} M-G^T) \beta_{1}-(A_0^T-{\Lambda_1^0} M_0)\beta_0 \nonumber\\
\quad\quad -(G^T-\Lambda_aM)\beta_2+\Phi_{0}M_0\alpha_0^0 +\Phi_{b}
M\alpha_1-\Gamma_1^TQ\eta,\\
 \dot{\beta}_1=(\Phi_{1}M-A^T) \beta_{1}+\Phi_{a}^TM_0\alpha_0^0 +
 \Phi_{2}M\alpha_1+Q\eta, \\
 \dot{\beta}_2=(\Phi_{2}^T M-F^T) \beta_{1}-(F_0^T-{\Lambda_2^0}^T M_0)\beta_0 \nonumber\\
\quad\quad -(A^T+F^T-(\Lambda_1+\Lambda_2^T)M)\beta_2 \nonumber\\
\quad\quad+\Phi_{b}^TM_0\alpha_0^0+\Phi_{3}M\alpha_1
-\Gamma_2^TQ\eta,
\end{array}\right.
\end{align}
where
\begin{align}
& \beta_0(T)=\Gamma_{1f}^TQ_f\eta_f ,\
 \beta_1(T)=-Q_f\eta_f ,\
 \beta_2(T)=\Gamma_{2f}^TQ_f\eta_f.  \nonumber
\end{align}
The optimal control law is given by
$$
u_1^*(t)=-R^{-1}B^T (\Phi_a^T X_0 +\Phi_1X_1+\Phi_2 \overline X +\beta_1).
$$

\begin{theorem} \label{theorem:PP}
We have
\begin{align}
&\mathbb P_0=\left[\begin{array}{ll}
\Lambda_{1}^0&\Lambda_{2}^0\\
{\Lambda_2^0}^T&\Lambda_{3}^0
\end{array}\right],\quad
\mathbb S_0=\left[\begin{array}{lll}
\alpha_0^0\\
\alpha_1^0
\end{array}\right],\nonumber
\end{align}
and \begin{align}
&\mathbb P =\left[\begin{array}{ccc}
\Lambda_{0},&\Lambda_{a},&\Lambda_{b}\\
{\Lambda_{a}^T},&\Lambda_{1},&\Lambda_{2}\\
{\Lambda_{b}^T},&\Lambda_{2}^T,&\Lambda_{3}
\end{array}\right],\quad
\mathbb S=\left[\begin{array}{cc}
\alpha_0\\
\alpha_1\\
\alpha_2
\end{array}\right].\nonumber
\end{align}
\end{theorem}

\proof
 We can directly show that
$$(\Phi_1^0,\Phi_2^0,\Phi_3^0,\Phi_0,\Phi_1,\Phi_2,\Phi_3,\Phi_a,\Phi_b)$$ is a solution of \eqref{eqn1*}. Similarly, $(\beta_0^0,\beta_1^0,\beta_0,\beta_1,\beta_2)$ satisfies the ODE of $(\alpha_0^0,\alpha_1^0,\alpha_0,\alpha_1,\alpha_2)$. Therefore, we obtain the representation of $(\mathbb P_0,\mathbb P,\mathbb S_0,\mathbb S)$.\hfill$\Box$

By Theorem \ref{theorem:PP},
after appropriate arrangement
the matrix functions in the solution of \eqref{eqn1*}
have the interpretation as the solutions of two Riccati-like equations.

It is now clear that $(u_0^*,u_1^*)$ agrees with
$(\check u_0, \check u_i)$ given by \eqref{u0}-\eqref{ui}.
Then we have the following interpretation on $(\check u_0, \check u_i)$
within an infinite population of minor players. First, $\check u_0$ is the best response with respect to $\overline X$; second, $\check u_i$ is the best response with respect to $({\overline X}, x_0, \check u_0)$; and finally, $\overline X$ is generated by the infinite number of minor players applying their best responses.
This suggests a consistent mean field approximation, which is well known in the fixed point approach of
mean field games (Caines, Huang, and Malham\'e, 2017).
In our mean field limit here, the consistent mean field approximation is a derived property. This is in contrast to the major player model in (Huang, 2010; Carmona and Zhu, 2016), where consistent mean
field approximations are imposed as
a requirement at the beginning so that individual strategies can be determined.

\section{Numerical examples}
\label{sec:num}
We have seen that testing asymptotic solvability reduces to checking the solution of \eqref{eqn1*} on $[0,T]$. It is generally infeasible to solve \eqref{eqn1*} analytically. Its numerical solution provides a practical means to check asymptotic solvability.
 \begin{figure}[t]
\begin{center}
\begin{tabular}{c}
\psfig{file=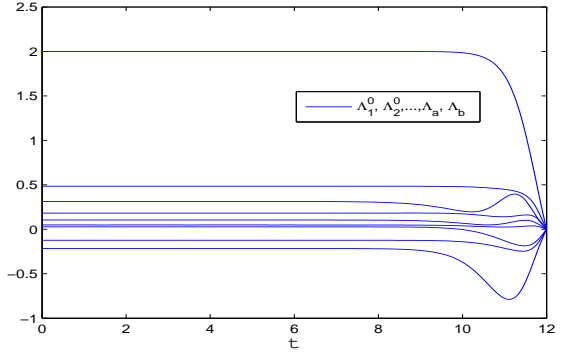, width=3.4in, height=2in}
\end{tabular}
\end{center}
\caption{The numerical  solution  of \eqref{eqn1*} in Example
\ref{example:Ex1}  } \label{fig1}
\end{figure}

\begin{figure}[t]
\begin{center}
\begin{tabular}{c}
\psfig{file=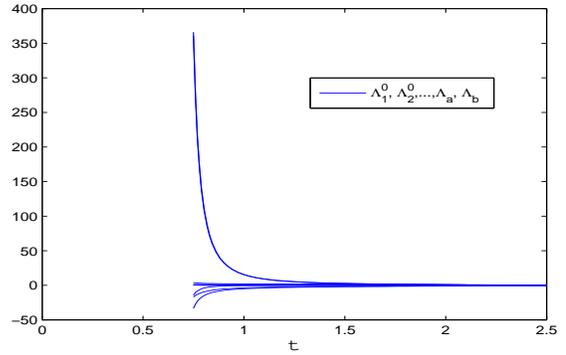, width=3.4in, height=2in}
\end{tabular}
\end{center}
\caption{The numerical  solution  of \eqref{eqn1*} in Example
\ref{example:Ex2}  } \label{fig2}
\end{figure}

\begin{example}\label{example:Ex1}
 The  parameters in \eqref{stateX0}-\eqref{costJi} are given by $A_0=1$, $ B_0=2$, $F_0=0.5$, $A=0.5$, $B=1$,  $F=0.2$, $ G=0.4$,
$Q_0=1$, $R_0=0.5$, $Q=2$, $ R=1$, $ \Gamma_0=0.8$, $\Gamma_1=0.3$, $\Gamma_2=0.5$, $Q_{0f}=Q_f=0$,  and $T=12$.
 We use ode45 of MatLab to numerically solve \eqref{eqn1*}  on $[0,T]$; see Fig. \ref{fig1}.
The existence of the solution suggests asymptotic solvability holds.
\end{example}

\begin{example} \label{example:Ex2}
$A_0=0.3$, $B_0=1$, $F_0=0.2$, $A=0.2$, $B=1$,  $F=1$, $G=-0.2$,
$Q_0=2$, $R_0=1$, $Q=1$, $R=1$, $\Gamma_0=0.8$, $\Gamma_1=0.1$, $\Gamma_2=1.2$, $Q_{0f}=Q_f=0$, and $T=2.5$.
The numerical solution of \eqref{eqn1*} has a finite escape time between $0.5$ and $1$ as shown in Fig. \ref{fig2}, suggesting no asymptotic
solvability.
\end{example}

\section{Concluding remarks}
\label{sec:con}

We study an asymptotic solvability problem for LQ Nash games involving a major player and $N$ minor players, where $N$ tends to infinity. We obtain the necessary and sufficient condition of asymptotic solvability via a system of Riccati ODEs and evaluate the equilibrium costs.
The system of Riccati ODEs has close relation with a limiting control model
of two players: the major player and a representative minor player.
For future work, it is of interest to generalize our analysis to deal  with leadership (Bensoussan et al, 2017), noisy measurements (Firoozi and Caines, 2015), and control constraints (Hu, Huang and Li, 2018).


\section*{Appendix A: Proof of Theorem \ref{theorem:Prep}}
\renewcommand{\theequation}{A.\arabic{equation}}
\setcounter{equation}{0}
\renewcommand{\thetheorem}{A.\arabic{theorem}}
\setcounter{theorem}{0}

\begin{lemma}
\label{lemma:Prep}
 Assume that \eqref{DE3_P0} and \eqref{DE3_P} have a solution
 $(\mbP_0(t), \cdots,
\mbP_N(t))$  on $[0,T]$. Then the following holds.

i) ${\mbP}_0(t)$  has the  representation
\begin{align}
{\mbP}_0(t)=\begin{bmatrix}
\Pi_1^0 &\Pi_2^0 &\Pi_2^0&\cdots &\Pi_2^0 \\
{\Pi_2^0}^T &\Pi_3^0 &\Pi_4^0&\cdots &\Pi_4^0\\
{\Pi_2^0}^T &\Pi_4^0&\Pi_3^0 &\cdots &\Pi_4^0\\
\vdots      & \vdots      &  \vdots     &\ddots    &\vdots \\
{\Pi_2^0}^T &\Pi_4^0 &\Pi_4^0&\cdots &\Pi_3^0
\end{bmatrix},\nonumber
\end{align}
where $\Pi_1^0(t)$, $\Pi_3^0(t)$ and $\Pi_4^0(t)$ are $n\times n$
symmetric matrix functions. The matrix $\Pi_4^0$ appears for $N^2-N$ times.

ii) ${\mbP}_1(t)$  has the  representation
\begin{align}
{\mbP}_1(t)=\begin{bmatrix}
\Pi_0   &\Pi_a   &\Pi_b  &\cdots   &\Pi_b \\
\Pi_a^T &\Pi_1   &\Pi_2 &\cdots   &\Pi_2\\
\Pi_b^T &\Pi_2^T &\Pi_3  &\cdots   &\Pi_4\\
\vdots     & \vdots    &  \vdots  &\ddots   &\vdots \\
\Pi_b^T &\Pi_2^T &\Pi_4  &\cdots   &\Pi_3
\end{bmatrix}\nonumber
\end{align}
where $\Pi_0(t)$, $\Pi_1(t)$, $\Pi_3(t)$ and $\Pi_4(t)$ are $n\times n$
symmetric matrix functions. The matrix $\Pi_4$ appears for $(N-1)(N-2)$ times.

iii) For $i>1$,  $\mbP_i(t)= J_{2,i+1}^T \mbP_1(t) J_{2,i+1}$.
\end{lemma}

\proof
i)  For $0\le l\le N$,  denote
$\mbP_l= (\mbP_l^{jk})_{1\le j,k\le N+1}, $
where $\mbP_l^{jk}$ is an $n\times n$ matrix.
Let $\hat{\mbP_l}=J_{23}^T \mbP_l J_{23}$.
By the method in (Huang and Zhou, 2018b, Lemma A.1) and elementary matrix computations, we can verify that
$(\hat{\mbP_0},\hat{\mbP_2},\hat{\mbP_1},\hat{\mbP_3},\cdots, \hat{\mbP_N})$ satisfies \eqref{DE3_P0} and \eqref{DE3_P}. Hence,$$\hat{\mbP_0}=\mbP_0,~\hat{\mbP_2}=\mbP_1,~\hat{\mbP_1}=\mbP_2,~\hat{\mbP_k}=\mbP_k,~k\geq3.$$
Then ${\mbP_2}=J_{23}^T \mbP_1 J_{23}$ and ${\mbP_0}=J_{23}^T \mbP_0 J_{23}$, and we obtain $\mbP_0^{22}=\mbP_0^{33},~\mbP_0^{12}=\mbP_0^{13}$.

Taking $J_{k,k+1},~k\geq3$ in place of $J_{23}$ and following the method in (Huang and Zhou, 2018b), we obtain the representation of $\mbP_0$.

ii) Now denote $\hat{\mbP_l}=J_{34}^T \mbP_l J_{34}$, and we can verify that
$$(\hat{\mbP_0},\hat{\mbP_1},\hat{\mbP_3},\hat{\mbP_2},\hat{\mbP_4},\cdots, \hat{\mbP_N})$$
is a solution of \eqref{DE3_P0} and \eqref{DE3_P}. Hence,
$$\hat{\mbP_1}=\mbP_1,~\hat{\mbP_3}=\mbP_2,~\hat{\mbP_2}=\mbP_3.$$
This yields $\mbP_1^{13}=\mbP_1^{14}$, $\mbP_1^{23}=\mbP_1^{24}$ and $\mbP_1^{33}=\mbP_1^{44}$.
In addition, $\mbP_1^{43}=\mbP_1^{34}$. Since $\mbP_1$ is symmetric, $(\mbP_1^{43})^T=\mbP_1^{34}$. So $\mbP_1^{34}$ is symmetric. Similarly, by the relation $J_{45}^T\mbP_1J_{45}=\mbP_1,$ we obtain $\mbP_1^{34}=\mbP_1^{35}$. Now, repeatedly using the relation
$J_{k,k+1}^T\mbP_1J_{k,k+1}=\mbP_1$ for all $k\ge 3$,
 we obtain the representation of $\mbP_1$. Note that $\Pi_4$ is symmetric.

iii) This equality can be shown as in the case $i=2$ in the proof of part i).\hfill$\Box$

{\bf Proof of Theorem \ref{theorem:Prep}:}
By Lemma \ref{lemma:Prep}, we have
\begin{align}
&\dot{\Pi_1^0}(t)  =   \Pi_1^0  M_0\Pi_1^0 +N\Pi_2^0 M\Pi_a^T
                          +N\Pi_a  M{\Pi_2^0}^T  \nonumber \\
                          &\quad\qquad - (\Pi_1^0 A_0+A_0^T\Pi_1^0 ) \nonumber \\
                          &\quad\qquad -N(\Pi_2^0 G+ G^T {\Pi_2^0}^T )-Q_{0},
                             \nonumber \\
                  &\Pi_1^0(T)=Q_{0f},\label{Pi10}
\end{align}
and
\begin{align}\label{Pi20}
&\dot{\Pi_2^0}(t) =  \Pi_1^0 M_0\Pi_2^0 +\Pi_2^0  M\Pi_1 +\Pi_a M\Pi_3^0 \nonumber \\
             &\quad\qquad   +(N-1)(\Pi_2^0 M\Pi_2  + \Pi_a M\Pi_4^0)  \nonumber \\
             &\quad\qquad   -(\Pi_1^0  \tfrac{F_0}{N}+\Pi_2^0 F+ \Pi_2^0 A) \nonumber \\
             &\quad\qquad  -(A_0^T\Pi_2^0 +G^T\Pi_3^0 +(N-1)G^T\Pi_4^0 )+ Q_{0}\tfrac{\Gamma_{0}}{N},\nonumber \\
&\Pi_2^0(T)=-Q_{0f}\tfrac{\Gamma_{0f}}{N},
\end{align}
and
\begin{align}\label{Pi30}
&\dot{\Pi_3^0}(t) = {\Pi_2^0}^T M_0\Pi_2^0 +{\Pi_3^0} M\Pi_1 +\Pi_1 M\Pi_3^0  \nonumber \\
                       &\quad\qquad +(N-1)(\Pi_4^0 M\Pi_2 +\Pi_2^T M\Pi_4^0)  \nonumber\\
                       &\quad\qquad - \tfrac{1}{N}({\Pi_2^0}^T F_0+F_0^T\Pi_2^0)  \nonumber\\
                       &\quad\qquad -(1-\tfrac{1}{N}) ( \Pi_4^0 {F}+ F^T\Pi_4^0 ) \nonumber\\
                       &\quad\qquad - \Big( \Pi_3^0 (A+\tfrac{F}{N})+ (A^T+\tfrac{F^T}{N})\Pi_3^0 \Big) - \tfrac{\Gamma_{0}^T}{N}Q_{0}\tfrac{\Gamma_{0}}{N},\nonumber \\
& \Pi_3^0(T)=\tfrac{\Gamma_{0f}^T}{N}Q_{0f}\tfrac{\Gamma_{0f}}{N},
\end{align}
and
\begin{align}\label{Pi40}
&\dot{\Pi_4^0}(t) ={\Pi_2^0}^T M_0\Pi_2^0 +{\Pi_4^0} M\Pi_1 +\Pi_1 M{\Pi_4^0} \nonumber \\
               &\quad\qquad + {\Pi_3^0} M\Pi_2 +\Pi_2^T M{\Pi_3^0}  \nonumber \\
               &\quad\qquad + (N-2)({\Pi_4^0} M\Pi_2 + \Pi_2^T M{\Pi_4^0})  \nonumber \\
               &\quad\qquad -\tfrac{1}{N} ({\Pi_2^0}^T{F_0}+F_0^T{\Pi_2^0} )
               -\tfrac{1}{N} ( \Pi_3^0 {F} + {F^T} \Pi_3^0 )  \nonumber \\
               &\quad\qquad - \Big(\Pi_4^0 (A+\tfrac{N-1}{N}F) + (A^T+\tfrac{N-1}{N}F^T)\Pi_4^0  \Big)  \nonumber \\
               &\quad\qquad - \tfrac{\Gamma_{0}^T}{N}Q_{0}\tfrac{\Gamma_{0}}{N},   \nonumber \\
& \Pi_4^0(T)=\tfrac{\Gamma_{0f}^T}{N}Q_{0f}\tfrac{\Gamma_{0f}}{N}.
\end{align}
We have $\Pi_3^0(T)=\Pi_4^0(T)$, and
\begin{align}
\dot{\Pi_3^0}(t)- \dot{\Pi_4^0}(t)=&(\Pi_3^0 - \Pi_4^0 )(M\Pi_1 -M\Pi_2 -A)\nonumber\\
&+(\Pi_1M-\Pi_2^T M-A^T)(\Pi_3^0 - \Pi_4^0 ).\nonumber
\end{align}
It follows that
\begin{align}   \label{Pi30=Pi40}
\Pi_3^0(t)=\Pi_4^0(t)\quad  \forall t\in [0,T].
\end{align}

By Lemma \ref{lemma:Prep}, we have
\begin{align}\label{Pi0}
&\dot{\Pi_0}(t)  =  \Pi_aM\Pi_a^T +\Pi_0  M_0\Pi_1^0  + \Pi_1^0  M_0 \Pi_0  \nonumber \\
               &\quad\qquad  - (\Pi_0 A_0 + A_0^T\Pi_0 )  - (\Pi_a G+ G^T \Pi_a^T )\nonumber\\
               &\quad\qquad  -(N-1) (\Pi_b G+G^T \Pi_b^T )\nonumber\\
               &\quad\qquad  + (N-1)(\Pi_b M\Pi_a^T
               + \Pi_a M\Pi_b^T) \nonumber \\
               &\quad\qquad  - \Gamma_1^T Q \Gamma_1,\nonumber\\
                  &\Pi_0(T)=\Gamma_{1f}^TQ_{f}\Gamma_{1f},
\end{align}
and  
\begin{align}\label{Pi1}
&\dot{\Pi_1}(t) = \Pi_1 M\Pi_1 + {\Pi_a^T}M_0\Pi_2^0 +{\Pi_2^0}^TM_0\Pi_a  \nonumber \\
                &\quad\qquad  + (N-1)(\Pi_2 M\Pi_2 +\Pi_2^T M\Pi_2^T)  \nonumber\\
                &\quad\qquad  - \tfrac{1}{N}({\Pi_a^T} {F_0}+{F_0^T}\Pi_a)  \nonumber\\
                &\quad\qquad  - \Big({\Pi_1} (A+\tfrac{F}{N})+(A^T+
                \tfrac{F^T}{N})\Pi_1 \Big)  \nonumber\\
                &\quad\qquad - (1-\tfrac{1}{N}) ( \Pi_2 {F}+ {F^T}\Pi_2^T) \nonumber\\
                &\quad\qquad  - (I- \tfrac{\Gamma_{2}^T}{N})Q(I-\tfrac{\Gamma_2}{N}),  \nonumber \\
& \Pi_1(T)=(I- \tfrac{\Gamma_{2f}^T}{N})Q_{f}(I-\tfrac{\Gamma_{2f}}{N}),
\end{align}
and
\begin{align}\label{Pi2}
&\dot{\Pi_2}(t) ={\Pi_a^T}M_0\Pi_2^0+{\Pi_2^0}^TM_0{\Pi_b}\nonumber \\
               &\quad\qquad + {\Pi_1 }M\Pi_2 +{\Pi_2} M\Pi_1 +\Pi_2^T M{\Pi_3} \nonumber \\
               &\quad\qquad + (N-2)({\Pi_2} M\Pi_2 + \Pi_2^T M{\Pi_4})  \nonumber \\
               &\quad\qquad -\tfrac{1}{N} (\Pi_a^T {F_0}+{F_0^T}{\Pi_b})  \nonumber \\
               &\quad\qquad - \Big(\Pi_2 (A+\tfrac{N-1}{N}F) + (A^T+\tfrac{F^T}{N})\Pi_2  \Big)  \nonumber \\
               &\quad\qquad -  (1-\tfrac{2}{N}){F^T}\Pi_4-\tfrac{1}{N} ( \Pi_1 {F} + {F^T} \Pi_3  ) \nonumber \\
               &\quad\qquad + (I- \tfrac{\Gamma_{2}^T}{N})Q\tfrac{\Gamma_2}{N}   ,\nonumber \\
& \Pi_2(T)=-(I- \tfrac{\Gamma_{2f}^T}{N})Q_{f}\tfrac{\Gamma_{2f}}{N},
\end{align}
and
\begin{align}\label{Pi3}
&\dot{\Pi_3}(t) ={\Pi_b^T} M_0\Pi_2^0 +{\Pi_2^0}^T M_0{\Pi_b}   \nonumber \\
               &\quad\qquad + \Pi_2^T M\Pi_2 +\Pi_3 M\Pi_1 +\Pi_1 M\Pi_3  \nonumber \\
               &\quad\qquad + (N-2)({\Pi_4} M\Pi_2 + \Pi_2^T M{\Pi_4})  \nonumber \\
               &\quad\qquad -\tfrac{1}{N} (\Pi_b^T {F_0}+{F_0^T}{\Pi_b} + \Pi_2^T {F} + {F^T} \Pi_2  ) \nonumber \\
               &\quad\qquad - \Big(\Pi_3 (A+\tfrac{F}{N}) + (A^T+\tfrac{F^T}{N})\Pi_3  \Big)  \nonumber \\
               &\quad\qquad -(1-\tfrac{2}{N}) (\Pi_4 {F}+ {F^T}\Pi_4 ) -  \tfrac{\Gamma_{2}^T}{N}Q\tfrac{\Gamma_2}{N},   \nonumber \\
& \Pi_3(T)=\tfrac{\Gamma_{2f}^T}{N}Q_{f}\tfrac{\Gamma_{2f}}{N} ,
\end{align}
and
\begin{align}\label{Pi4}
&\dot{\Pi_4}(t) ={\Pi_b^T} M_0\Pi_2^0 +{\Pi_2^0}^T M_0{\Pi_b}+ {\Pi_2^T }M\Pi_2 +{\Pi_3} M\Pi_2 \nonumber \\
               &\quad\qquad +\Pi_2^T M{\Pi_3} + {\Pi_4 }M\Pi_1 +{\Pi_1} M\Pi_4   \nonumber \\
               &\quad\qquad + (N-3)({\Pi_4} M\Pi_2 + \Pi_2^T M{\Pi_4})  \nonumber \\
               &\quad\qquad -\tfrac{1}{N} (\Pi_b^T {F_0}+{F_0^T}{\Pi_b})  \nonumber \\
               &\quad\qquad - \Big(\Pi_4 (A+\tfrac{N-2}{N}F) + (A^T+\tfrac{N-2}{N}F^T)\Pi_4  \Big)  \nonumber \\
               &\quad\qquad - \tfrac{1}{N}( \Pi_2^T {F} + {F^T}
               \Pi_2  + \Pi_3{F} + {F^T} \Pi_3  )
               - \tfrac{\Gamma_{2}^T}{N}Q\tfrac{\Gamma_2}{N} ,  \nonumber \\
& \Pi_4(T)=\tfrac{\Gamma_{2f}^T}{N}Q_{f}\tfrac{\Gamma_{2f}}{N},
\end{align}
and
\begin{align}\label{Pia}
&\dot{\Pi_a}(t) =  \Pi_0  M_0 \Pi_2^0 + \Pi_1^0 M_0\Pi_a  + \Pi_a M\Pi_1 \nonumber \\
             &\quad\qquad   +(N-1)(\Pi_b M\Pi_2  + \Pi_a M\Pi_2^T)  \nonumber \\
             &\quad\qquad   -\Big(\Pi_0  \tfrac{F_0}{N} + \Pi_a (A+\tfrac{F}{N})+ (N-1)\Pi_b \tfrac{F}{N}\Big) \nonumber \\
             &\quad\qquad   -\big(A_0^T\Pi_a +G^T\Pi_1 +(N-1)G^T\Pi_2^T\big)\nonumber \\
             &\quad\qquad   +\Gamma_1^TQ(I-\tfrac{\Gamma_2}{N}),\nonumber \\
&\Pi_a(T)=-\Gamma_{1f}^TQ_{f}(I-\tfrac{\Gamma_{2f}}{N}),
\end{align}
and
\begin{align}\label{Pib}
&\dot{\Pi_b}(t) = {\Pi_0} M_0\Pi_2^0 +{\Pi_1^0} M_0{\Pi_b} +{\Pi_b} M{\Pi_1}   \nonumber \\
               &\quad\qquad + \Pi_a M\Pi_2 +\Pi_a M\Pi_3  \nonumber \\
               &\quad\qquad + (N-2)({\Pi_b} M\Pi_2 + \Pi_a M{\Pi_4})  \nonumber \\
               &\quad\qquad-\Big({\Pi_0} \tfrac{F_0}{N}+\Pi_a \tfrac{F}{N}+{\Pi_b} A+(N-1){\Pi_b} \tfrac{F}{N}\Big)  \nonumber \\
               &\quad\qquad - ( A_0^T\Pi_b + G^T \Pi_2 + G^T \Pi_3 + (N-2)G^T \Pi_4 ) \nonumber \\
               &\quad\qquad -  \Gamma_{1}^TQ\tfrac{\Gamma_2}{N},   \nonumber \\
& \Pi_b(T)=\Gamma_{1f}^TQ_f\tfrac{\Gamma_{2f}}{N} .
\end{align}

We have  $\Pi_3(T)=\Pi_4(T)$, and
\begin{align}
\dot{\Pi_3}(t)- \dot{\Pi_4}(t)=&(\Pi_3 - \Pi_4 )(M\Pi_1 -M\Pi_2 -A)\nonumber\\
&+(\Pi_1M-\Pi_2^T M-A^T)(\Pi_3 -\Pi_4 ).\nonumber
\end{align}
Therefore,
$\Pi_3(t)=\Pi_4(t)$ for all  $ t\in [0,T]$.
This completes the proof.
\endproof

\section*{Appendix B: Proof of Theorem \ref{theorem:AS}}
\renewcommand{\theequation}{B.\arabic{equation}}
\setcounter{equation}{0}
\renewcommand{\thetheorem}{B.\arabic{theorem}}
\setcounter{theorem}{0}

Step 1.
By \eqref{Pi10}, \eqref{Pi30} and \eqref{Pi30=Pi40}, we determine
\begin{align}
&\dot{\Lambda}_1^{0N} =  \Lambda_1^{0N}M_0 \Lambda_1^{0N}+\Lambda_2^{0N}M{\Lambda_a^N}^T+\Lambda_a^N M (\Lambda_2^{0N})^T \nonumber\\
&\quad\quad\quad-(\Lambda_1^{0N}A_0+A_0^T\Lambda_1^{0N})\nonumber\\
&\quad\quad \quad-(\Lambda_2^{0N}G+G^T(\Lambda_2^{0N})^T)-Q_0 ,\label{od10}
\end{align}
\begin{align}
&\dot{\Lambda}_2^{0N} = \Lambda_1^{0N}M_0 \Lambda_2^{0N}+\Lambda_2^{0N}(M\Lambda_1^{N}+M \Lambda_2^{N}-F-A) \nonumber\\
&\quad\quad\quad-
A_0^T \Lambda_2^{0N}-\Lambda_1^{0N}F_0 +(\Lambda_a^{N}M-G^T)\Lambda_3^{0N}\nonumber\\
&\quad\quad\quad +Q_0\Gamma_0+g_2^0(1/N,\Lambda_2^{0N},\Lambda_2^{N}), \label{d11}
\end{align}
\begin{align}
&\dot{\Lambda}_3^{0N} = (\Lambda_2^{0N})^TM_0 \Lambda_2^{0N}+\Lambda_3^{0N}M\Lambda_1^{N}+\Lambda_1^{N} M \Lambda_3^{0N}- (\Lambda_2^{0N})^TF_0\nonumber\\
&\quad\quad\quad-F_0^T\Lambda_2^{0N}-\Lambda_3^{0N}(A+F)- (A^T+F^T)\Lambda_3^{0N} \nonumber\\
&\quad\quad\quad +\Lambda_3^{0N}M \Lambda_2^{N}+{\Lambda_2^{N}}^TM \Lambda_3^{0N}\nonumber\\
&\quad\quad\quad -\Gamma_0^TQ_0\Gamma_0+g_3^0(1/N,
\Lambda_3^{0N},\Lambda_2^{N}), \label{d1130}
\end{align}
where $\Lambda_1^{0N}(T)=Q_{0f}$, $\Lambda_2^{0N}(T)=-Q_{0f}\Gamma_{0f}$,
$\Lambda_3^{0N}(T)=\Gamma_{0f}^TQ_{0f}\Gamma_{0f}$, and
\begin{align}
&g_2^0(1/N,\Lambda_2^{0N},\Lambda_2^{N}
)=-(1/N)(\Lambda_2^{0N}M \Lambda_2^{N}).\nonumber
\end{align}
For reasons of space, the expression of $g_3^0$ is not displayed.
We further obtain
\begin{align}
&\dot{\Lambda}_0^{N} = \Lambda_a^{N}M{\Lambda_a^{N}}^T+\Lambda_0^{N}M_0\Lambda_1^{0N}+\Lambda_1^{0N} M_0 \Lambda_0^{N}\nonumber\\
&\quad\quad\quad -\Lambda_0^{N}A_0-A_0^T\Lambda_0^{N}-(\Lambda_a^{N}+\Lambda_b^{N})G- G^T(\Lambda_a^{N}+ \Lambda_b^{N} )^T\nonumber\\
&\quad\quad\quad +\Lambda_b^{N}M {\Lambda_a^{N}}^T+\Lambda_a^{N}M{\Lambda_b^{N}}^T\nonumber\\
&\quad\quad\quad -\Gamma_1^TQ\Gamma_1+g_0(1/N,\Lambda_a^{N},
\Lambda_b^{N}), \label{d110}\\
&\Lambda_0^N(T)=\Gamma_{1f}^TQ_f\Gamma_{1f},\nonumber  
\end{align}
\begin{align}
&\dot{\Lambda}_1^{N} = \Lambda_1^{N}M\Lambda_1^{N}-\Lambda_1^{N}A-A^T\Lambda_1^{N}-Q\nonumber\\
&\quad\quad\quad +g_1(1/N,\Lambda_2^{0N},\Lambda_1^{N},\Lambda_2^{N},\Lambda_a^{N}),    \label{d111}\\
&\Lambda_1^N(T)=(I- \tfrac{\Gamma_{2f}^T}{N})Q_{f}(I-\tfrac{\Gamma_{2f}}{N}),\nonumber
\end{align}
\begin{align}
&\dot{\Lambda}_2^{N} = {\Lambda_a^{N}}^T(M_0\Lambda_2^{0N}-F_0)+\Lambda_1^{N}M\Lambda_2^{N}
+\Lambda_2^{N}
M\Lambda_1^{N}\nonumber\\
&\quad\quad\quad-\Lambda_2^{N}(A+F)-A^T\Lambda_2^{N}-\Lambda_1^{N}F
+\Lambda_2^{N}M\Lambda_2^{N}\nonumber\\
&\quad\quad\quad+Q\Gamma_2 +g_2(1/N,\Lambda_2^{0N},\Lambda_2^{N},\Lambda_3^{N},\Lambda_b^{N}),    \label{d112}\\
&\Lambda_2^N(T)=-(I-\tfrac{\Gamma_{2f}^T}{N})Q_f\Gamma_{2f},\nonumber
\end{align}
\begin{align}
&\dot{\Lambda}_3^{N} = {\Lambda_b^{N}}^TM_0\Lambda_2^{0N}+(\Lambda_2^{0N})^TM_0\Lambda_b^{N}+{\Lambda_2^{N}}^TM\Lambda_2^{N}\nonumber\\
&\quad\quad\quad +\Lambda_3^{N}M\Lambda_1^{N}+\Lambda_1^{N}M\Lambda_3^{N} + (\Lambda_3^{N}M\Lambda_2^{N}+{\Lambda_2^{N}}^TM\Lambda_3^{N}
) \nonumber\\
&\quad\quad\quad-{\Lambda_b^{N}}^TF_0-F_0^T\Lambda_b^{N} -{\Lambda_2^{N}}^TF-F^T\Lambda_2^{N}  \nonumber\\
&\quad\quad\quad -\Lambda_3^{N}(A+F)-(A^T+F^T)\Lambda_3^{N}- \Gamma_2^T Q\Gamma_2 \nonumber\\
&\quad\quad\quad  +g_3(1/N,\Lambda_2^{N},\Lambda_3^{N}),    \label{d3n}\\
&\Lambda_3^N(T)=\Gamma_{2f}^T Q_f\Gamma_{2f},\nonumber
\end{align}
\begin{align}
&\dot{\Lambda}_a^{N} = \Lambda_1^{0N}M_0\Lambda_a^{N}
+\Lambda_a^{N}M\Lambda_1^{N}+\Lambda_a^NM{\Lambda_2^N}^T\nonumber\\
&\quad\quad\quad-\Lambda_a^{N}A -A_0^T\Lambda_a^{N}-G^T(\Lambda_1^{N}+{\Lambda_2^N}^T)
   \nonumber\\
&\quad\quad\quad +\Gamma_1^TQ+g_a(1/N,\Lambda_2^{0N},\Lambda_0^{N},\Lambda_2^{N},
\Lambda_a^{N},\Lambda_b^{N}) ,   \label{dan}\\
&\Lambda_a^N(T)=-\Gamma_{1f}^T Q_f(I-\tfrac{\Gamma_{2f}}{N}),\nonumber
\end{align}
\begin{align}
&\dot{\Lambda}_b^{N} = \Lambda_0^{N}M_0\Lambda_2^{0N}+\Lambda_1^{0N}M_0\Lambda_b^{N}+\Lambda_b^{N}M\Lambda_1^{N}+\Lambda_a^{N}M\Lambda_2^{N}\nonumber\\
&\quad\quad\quad+\Lambda_b^{N}M\Lambda_2^{N}+\Lambda_a^{N}M\Lambda_3^{N}-\Lambda_0^NF_0-\Lambda_a^NF\nonumber\\
&\quad\quad\quad-\Lambda_b^N(A+F)-A_0^T\Lambda_b^{N}-G^T(\Lambda_2^N+\Lambda_3^N) \nonumber\\
&\quad\quad\quad-\Gamma_1^T Q\Gamma_2+
g_b(1/N,\Lambda_2^{N},\Lambda_3^{N},\Lambda_a^{N},\Lambda_b^{N}),    \label{dbn}\\
&\Lambda_b^N(T)=\Gamma_{1f}^T Q_f\Gamma_{2f},
\nonumber\end{align}
where $g_0, \cdots, g_b$ are not displayed and are compactly of $O(1/N)$.


Step 2. Proof of Theorem \ref{theorem:AS}.

Denote $\xi^N=(\Lambda_1^{0N}, \Lambda_2^{0N}, \cdots, \Lambda_b^N )$ for \eqref{L0P}, and we view each of $g_2^0$, $g_3^0$, $g_0$, $g_1$, $g_2$, $g_3$,   $g_a$, $g_b$ as a function of $ \xi^N$ with parameter $1/N$.
 They are all compactly of $O(1/N)$.
For some $C>0$,
we further have
\begin{align}
&|(\Lambda_1^{0N}(T)-\Lambda_1^0(T), \Lambda_2^{0N}(T)-\Lambda_2^{0}(T), \ldots, \Lambda_b^N(T)-\Lambda_b(T))|\nonumber \\
&\le C/N. \nonumber
\end{align}
Subsequently, we view the ODE system \eqref{od10}-\eqref{dbn} as a slightly perturbed form of \eqref{eqn1*}. The remaining proof is similar to that of (Huang and Zhou, 2018b, Theorem 5) and we only give its sketch.
If asymptotic solvability holds, we solve \eqref{od10}-\eqref{dbn}
for all sufficiently large $N$. By taking
some increasing subsequence of population sizes
$N_1<N_2<\cdots$, we can ensure that as $k\to \infty$,
their solutions
$\{\xi^{N_k}(t), k=1, 2, \cdots\}$ have a
 limit
 as a vector function on $[0,T]$
which satisfies the limit ODE \eqref{eqn1*}  on $[0,T]$. Conversely, if \eqref{eqn1*} has a solution on $[0,T]$, there exists $N_0>0$ such that \eqref{od10}-\eqref{dbn} has a solution on $[0,T]$ for all  $N\ge N_0$; all these solutions are uniformly bounded. Accordingly we obtain $(\mbP_0, \ldots, \mbP_N)$ to satisfy \eqref{P0P1} for all $N\ge N_0$. So asymptotic solvability holds.~\endproof

\section*{Appendix C }
\renewcommand{\theequation}{C.\arabic{equation}}
\setcounter{equation}{0}
\renewcommand{\thetheorem}{C.\arabic{theorem}}
\setcounter{theorem}{0}

In view of \eqref{DE3_S0} and \eqref{DE3_S}, by Proposition \ref{prop:S} we have
\begin{align}
\left\{
\begin{array}{l}
\dot{\theta}_0^{0}(t) = -A_0^T\theta_0^0-NG^T\theta_1^0+
\Pi_1^0M_0\theta_0^0+N\Pi_aM\theta_1^0 \\
\quad\quad\quad +N\Pi_2^0M\theta_1+Q_0\eta_0,\\
\dot{\theta}_1^{0}(t) = -\frac{F_0^T}{N}\theta_0^{0}-A^T\theta_1^{0}-F^T\theta_1^{0}
+{\Pi_2^0}^TM_0\theta_0^{0}+\Pi_1M\theta_1^{0} \\
\quad\quad\quad+(N-1)\Pi_2^TM\theta_1^0+N\Pi_3^0M\theta_1-\frac{\Gamma_0^T}{N}Q_0\eta_0,
\end{array}
\right. \label{si0}
  \end{align}
where
$$
 \theta_0^{0}(T)=-Q_{0f}\eta_{0f},\qquad
 \theta_1^{0}(T)=\tfrac{\Gamma_{0f}^T}{N}Q_{0f}\eta_{0f},\nonumber
$$
and
\begin{align}
\left\{
\begin{array}{l}
\dot{\theta}_0(t) = -A_0^T\theta_0-G^T\theta_1-(N-1)G^T\theta_2+\Pi_1^0M_0\theta_0
 \\
\quad\qquad+\Pi_0M_0\theta_0^0 +\Pi_aM\theta_1+(N-1)\Pi_aM\theta_2 \\
\qquad\quad+(N-1)\Pi_bM\theta_1-\Gamma_1^TQ\eta,\\
\dot{\theta}_1(t) = -\frac{F_0^T}{N}\theta_0-(A^T+\frac{F^T}{N})\theta_1-(N-1)
\frac{F^T}{N}\theta_2+{\Pi_2^0}^TM_0\theta_0\\
\quad\qquad +\Pi_a^TM_0\theta_0^{0}+\Pi_1M\theta_1+(N-1)\Pi_2^TM\theta_2
\\
\quad\qquad+(N-1)
\Pi_2M\theta_1
 +(I-\frac{\Gamma_2^T}{N})Q\eta,\\
\dot{\theta}_2(t) = -\frac{F_0^T}{N}\theta_0-\frac{F^T}{N}\theta_1-A^T\theta_2-(N-1)
\frac{F^T}{N}\theta_2+{\Pi_2^0}^TM_0\theta_0\\
\quad\qquad +\Pi_b^TM_0\theta_0^{0}+\Pi_2^TM\theta_1+\Pi_1M\theta_2+(N-2)
\Pi_2^TM\theta_2\\
\quad\qquad +(N-1)\Pi_3M\theta_1-\frac{\Gamma_2^T}{N}Q\eta,
\end{array}
\right. \label{si12}
\end{align}
where
\begin{align}
\begin{cases}
 \theta_0(T)= \Gamma_{1f}^TQ_f\eta_f, \quad
 \theta_1(T)= -(I-\frac{\Gamma_{2f}^T}{N})Q_f\eta_f, \nonumber\\
 \theta_2(T)= \frac{\Gamma_{2f}^T}{N}Q_f\eta_f .\nonumber
\end{cases}
\end{align}

By \eqref{si0}, we have the relation
\begin{align}
&\dot{\alpha}_0^{0N} = (\Lambda_1^{0N}M_0-A_0^T)\alpha_0^{0N}
+(\Lambda_a^NM-G^T)\alpha_1^{0N}
 \nonumber\\
&\quad\quad +\Lambda_2^{0N}M\alpha_1^N+Q_0\eta_0,\nonumber \\
&\dot{\alpha}_1^{0N} =((\Lambda_2^{0N})^TM_0-F_0^T)\alpha_0^{0N}\nonumber\\
&\quad\quad+(\Lambda_1^N M+{\Lambda_2^N}^TM-A^T-F^T)\alpha_1^{0N} \nonumber\\
&\quad\quad
+\Lambda_3^{0N}M\alpha_1^N-\Gamma_0^TQ_0\eta_0 +h_1^0(1/N,\alpha_1^{0N},\Lambda_2^N),
\nonumber
\end{align}
where
 $\alpha _0^{0N}(T)=-Q_{0f}\eta_{0f}$,
 $\alpha _1^{0N}(T)=\Gamma_{0f}^TQ_{0f}\eta_{0f},$ 
and $ h_1^0(1/N,\alpha_1^{0N},\Lambda_2^N)=-\tfrac{1}{N}{\Lambda_2^N}^T\alpha_1^{0N}.$
By \eqref{si12}, we have
\begin{align*}
&\dot{\alpha}_0^{N} =(\Lambda_1^{0N}M_0 -A_0^T)\alpha_0^{N}+((\Lambda_a^N+\Lambda_b^N)M-G^T)\alpha_1^N \nonumber\\
&\quad\quad+(\Lambda_a^NM-G^T)\alpha_2^{N}
+\Lambda_0^NM_0\alpha_0^{0N} -\Gamma_1^TQ\eta\nonumber\\
&\quad\quad+h_0(1/N,\alpha_1^{N},\alpha_2^{N},\Lambda_a^{N},\Lambda_b^{N}),\\
& \dot{\alpha} _1^N = {\Lambda_a^N}^TM_0\alpha_0^{0N}-A^T\alpha_1^N+\Lambda_1^NM\alpha_1^N+\Lambda_2^NM\alpha_1^N+Q\eta\nonumber\\
&\quad\quad +h_1(\alpha_1^{0N},\alpha_0^N,\alpha_1^N,\alpha_2^N,
\Lambda_2^{0N},\Lambda_2^{N},\Lambda_a^{N}),\\
&\dot{\alpha}_2^N = ((\Lambda_2^{0N})^TM_0 -F_0^T)\alpha_0^{N}+
((\Lambda_3^N+{\Lambda_2^{N}}^T   ) M-F^T)\alpha_1^{N}\nonumber\\
&\quad\quad +((\Lambda_1^{N}+{\Lambda_2^{N}}^T)M-A^T-F^T)\alpha_2^{N}+
  {\Lambda_b^N}^TM_0 {\alpha_0^0}^N    \nonumber\\
&\quad\quad -\Gamma_2^TQ\eta+h_2(\alpha_1^N,\alpha_2^N,
\Lambda_2^{N},\Lambda_3^{N}),
\end{align*}
where the terminal condition is
\begin{align*}
 & \alpha _0^N(T)= \Gamma_{1f}^TQ_f\eta_f ,\quad
\alpha_1^N(T)= -(I-\tfrac{\Gamma_{2f}^T}{N})Q_f\eta_f, \\
 & \alpha _2^N(T)=\Gamma_{2f}^TQ_f\eta_f ,\nonumber
\end{align*} and $h_0$, $h_1$, $h_2$ are compactly of $O(1/N)$.


\end{document}